\documentclass{amsart}
\usepackage{Everything}

\title{Universal wide Aronszajn tree}
\author{Siiri Kivim\" aki}
\date{\today}
\thanks{This project has received funding from the European Research Council (ERC) under the
European Union’s Horizon 2020 research and innovation programme (grant agreement No
101020762) and Magnus Ehrnrooth Foundation.}

\let\oldtocsection=\tocsection

\let\oldtocsubsection=\tocsubsection

\let\oldtocsubsubsection=\tocsubsubsection

\renewcommand{\tocsection}[2]{\hspace{0em}\oldtocsection{#1}{#2}}
\renewcommand{\tocsubsection}[2]{\hspace{1.8em}\oldtocsubsection{#1}{#2}}
\renewcommand{\tocsubsubsection}[2]{\hspace{2em}\oldtocsubsubsection{#1}{#2}}

\makeatletter
\renewcommand\subsubsection{\@secnumfont}{\bfseries}
\renewcommand\subsubsection{\@startsection{subsubsection}{3}
  \z@{.5\linespacing\@plus.7\linespacing}{-.5em}%
  {\normalfont\bfseries}}

\tikzset{
  symbol/.style={
    draw=none,
    every to/.append style={
      edge node={node [sloped, allow upside down, auto=false]{$#1$}}}
  }
}

\begin{document}

\begin{abstract}
    A wide Aronszajn tree is a tree of size $\aleph_1$ with no uncountable branches. Assuming the consistency of the existence of a weakly compact cardinal, we show the consistency of the existence of a wide Aronszajn tree that is \textit{universal} in the sense that it contains an isomorphic copy of every wide Aronszajn tree. 
\end{abstract}

\maketitle

\tableofcontents

\section{Introduction}

This paper contributes to the study of universality problems. A universality problem asks, for a class $\KK$ equipped with a notion of embedding, does $\KK$ contain a \textit{universal object}? Such an object is a member of $\KK$ into which every other member embeds. For example, the rational line is a universal object in the class of countable linear orders, with respect to order-preserving embeddings, and the random graph is universal in the class of countable graphs, with respect to graph-homomorphisms. The picture changes radically when the objects are allowed to be uncountable. For instance, the existence of a universal model in an \textit{elementary class}, a class of the form $\Mod_\kappa(T)=\{M:M\models T\text{ and }|M|\leq\kappa\}$ equipped with elementary embeddings, is independent of the axioms of set theory whenever $T$ is unstable and $\kappa$ is an uncountable cardinal. 
This follows from the fact that, the existence of a saturated model of $T$ of size $\kappa$, which is always a universal object of $\Mod_\kappa(T)$, is equivalent to an instance of generalized continuum hypothesis, whenever $T$ is unstable (See Theorem 4.7 in \cite{shelah1990classification}). There is no such a combinatorial characterization for the existence of universal models and only little is known about them. 
Shelah showed in \cite{shelah1990universal} that there can exist a universal graph of size $\aleph_1$ even when the continuum hypothesis fails. This is an example of a scenario where a universal model of size $\aleph_1$ does exist, despite the fact that a saturated model of size $\aleph_1$ cannot not exist by the failure of the continuum hypothesis (to apply the above reasoning about elementry classes, in the case of graphs the relevant complete first-order theory is the theory of random graph). 
In \cite{mekler1990universal}, Mekler provides a general method for showing the consistency of the existence of a universal object of size $\aleph_1$ in a certain type of an amalgamation class together with the failure of the continuum hypothesis. In \cite{dvzamonja2006properties}, Shelah and Džamonja study properties of theories that would preclude the existence of a universal model whenever a saturated model does not exist.

In this paper we are interested in universal objects in non-elementary classes. We concentrate on a particular class of trees, the class of \textit{wide Aronszajn trees}. For basics about trees, we refer the reader to \cite{jech2003set}.

\begin{de}
    A \textbf{wide $\bs{\kappa}$-Aronszajn tree} is a tree of size $\kappa$ without linearly ordered subsets of size $\kappa$.
\end{de}

Wide $\kappa$-Aronszajn trees contrast with ordinary $\kappa$-Aronszajn trees, in which levels are required to be small. Contrary to ordinary $\kappa$-Aronszajn trees, wide $\kappa$-Aronszajn trees exist in every cardinality.

There are two relevant notions of embedding between trees:

\begin{de}
    Let $S$ and $T$ be trees. A function $f:S\to T$ is a
    \begin{enumerate}
        \item \textbf{weak embedding} if for any $s,t\in S$, 
        \[
        s<_St\implies f(s)<_Tf(t),
        \]
        \item \textbf{strong embedding} if it is level-preserving and for any $s,t\in S$, 
        \[
        s<_St\iff f(s)<_T f(t).
        \]
    \end{enumerate}
\end{de}

Weak embeddings need not be injective, but strong embeddings do. If there is a weak embedding $f$ from $S$ to $T$, then there is a weak embedding $\hat{f}$ from $S$ to $T$ that is also level-preserving. One simply redefines $\hat f(s)$ to be the unique predecessor of $f(s)$ at the height of $s$. 

It was shown by Kurepa \cite{kurepa1956ensembles} that the continuum hypothesis implies that there is no universal wide $\aleph_1$-Aronszajn tree. Later, Todorčević \cite{todorcevic2007lipschitz} and Džamonja and Shelah \cite{dvzamonja2021wide} showed that Martin's Axiom implies that there is no universal ordinary nor, respectively, wide $\aleph_1$-Aronszajn tree. All three results are with respect to weak embeddings.

In this paper it is shown that it is consistent that there is a wide $\aleph_1$-Aronszajn tree that is universal with respect to strong embeddings, assuming the consistency of a weakly compact cardinal. The proof is a forcing iteration of $\kappa$-strongly proper posets with side conditions. The construction builds on \cite{ben2023aronszajn}, where the consistency of the existence of a universal wide $\kappa^{++}$-Aronszajn tree is shown, for any cardinal $\kappa$. The result from \cite{ben2023aronszajn} leaves open the cases of $\aleph_1$ and successors of singulars. With this paper we complete the case of $\aleph_1$, leaving the successors of singulars still open. We conjecture that the result could be obtained without large cardinal hypotheses, by forcing over $L$ and replacing the use of the weakly compact cardinal by Devlin's version of $\Diamond$ from \cite{devlin1982combinatorial}.

\vv 

Our notation is standard and follows \cite{jech2003set}. A \textbf{tree} is a partial order that satisfies that the set of predecessors of any element is well-ordered. For a set $X$ and a cardinal $\kappa$, $\Pow_\kappa(X)$ is the set of subsets of $X$ of size $<\kappa$. The set $H_\kappa$ is the set of sets of hereditary cardinality $<\kappa$. If $\kappa$ is a regular cardinal, $(H_\kappa,\in)$ is a model of $\ZFC^-$ (all the axioms of set theory except maybe power set axiom). For a forcing poset $\P$ and $p,q\in\P$, we write $q\leq p$ to denote that the condition $q$ is stronger than $p$. By $\P/p$ we denote the set $\{q\in\P:q\leq p\}$. A model $M$ is \textbf{suitable for $\P$} if $M\elem H_\theta$ and $\P\in M$.

\vv 

\noindent\textbf{Aknowledgements.} This paper is part of my PhD thesis. I am grateful to professors Omer Ben-Neria and Menachem Magidor for various discussions, and to my supervisors professors Juliette Kennedy, Jouko Väänänen and Boban Veličković for guidance.

\vv

\subsection{Wide Aronszajn trees}\spa

\vv

The interest in wide Aronszajn trees grew from abstract model theory; a completely different area of logic than ordinary Aronszajn trees. We describe its starting point: Scott analysis for countable models.

\begin{de}
    Let $A$ and $B$ be models in a same signature and let $\beta$ be an ordinal. The game $\bs{\EF^\beta(A,B)}$ has two players and is played as follows: 
    A position in the game is a pair $(\alpha,\pi)$, where $\alpha\leq \beta$ is an ordinal and $\pi$ is a function from $A$ to $B$. The starting position is $(\beta,\emptyset)$. At position $(\alpha,\pi)$, player $\1$ chooses an ordinal $\alpha'$ and a point $a\in A$ (or a point $b\in B$), and player $\2$ responds by finding a point $b\in B$ (or a point $a\in A$). The next position is $(\alpha',\pi\cup\{(a,b)\})$. The game ends if either $\alpha'\not<\alpha$ in which case player $\1$ loses, or $\pi\cup\{(a,b)\}$ is not a partial isomorphism, in which case player $\2$ loses.
\end{de}

This variation of the Ehrenfeucht-Fraïssé game that uses a countable ordinal $\beta$ as a ``game clock" always results in a finite play, since there are no infinite descending sequences of ordinals. In particular $\EF^\beta$ is always determined. It is a powerful tool to analyse countable models.

The game $\EF^\beta$ is used to answer the question \begin{quote}
\textit{How different are non-isomorphic models?}
\end{quote}
Indeed, using the games $\EF^\beta$ it is possible to find a countable ordinal that can be seen as a parameter that measures the similarity of two countable models:
\begin{thm}[Scott \cite{scott2014logic}]\label{thm:scott1}
    Two countable models $A$ and $B$ are non-isomorphic if and only if there is a countable ordinal $\beta<\omega_1$ such that player $\2$ has a winning stategy in $\EF^\alpha(A,B)$ for $\alpha\leq\beta$ and player $\1$ has a winning strategy in $\EF^\gamma(A,B)$ for $\gamma>\beta$.
\end{thm}
It is said that $A$ and $B$ are ``isomorphic up to $\beta"$ if player $\2$ has a winning strategy in $\EF^\beta(A,B)$. Indeed, the games $\EF^\beta$ give rise to a ``grading" of the isomorphism relation; to $\aleph_1$ many weaker relations that approximate $\cong$ in the class of countable models.
Karp \cite{Karp1965-KARFQE} gives a logical characterization for these relations in terms of the infinitary logic $\LL_{\infty\omega}$ that extends the first-order logic by allowing infinitary conjunctions and disjunctions:

\begin{thm}[Karp \cite{Karp1965-KARFQE}]
    Let $A$ and $B$ be two structures in a countable relational signature. Player $\2$ has a winning strategy in $\EF^\beta(A,B)$ if and only if the models $A$ and $B$ are elementarily equivalent in the logic $\LL_{\infty\omega}$ up to quantifier rank $\beta$. 
\end{thm}

Thus if $A$ and $B$ are countable, then they are isomorphic if and only if they are elementarily equivalent in the logic $\LL_{\infty\omega}$. In particular, for non-isomorphic models there is always an $\LL_{\infty\omega}$-sentence (in fact, even an $\LL_{\omega_1\omega}$-sentence) that separates them.

The game $\EF^\beta$ is also used to find invariants and classify countable models. This application exploits the game-logic correspondence to its full power. Scott \cite{scott2014logic} realized that it is possible to code winning strategies in $\EF^\beta$ and for every countable model $A$, to write down an $\LL_{\omega_1\omega}$-sentence that characterizes $A$ up to isomorphism in the class of countable models:

\begin{thm}[Scott's isomorphism theorem \cite{scott2014logic}]\label{thm:scottiso} For any countable model $A$ there is an $\LL_{\omega_1\omega}$-sentence $\psi_A$ such that for any countable model $B$, $B$ is isomorphic to $A$ if and only if $B\models\psi_A$.
\end{thm}
The sentence $\psi_A$ is called the \textbf{Scott sentence} of $A$. 

It would be tempting to carry out the same analysis for uncountable models using the Ehrenfeucht-Fraïssé games corresponding to the logics $\LL_{\lambda\kappa}$, but this does not work. Indeed, there are linear orders of size $\aleph_1$ that are elementarily equivalent even in the stronger logic $\LL_{\infty\omega_1}$, and yet non-isomorphic. So, in general, it is impossible to separate non-isomorphic uncountable models using an $\LL_{\infty\omega_1}$-sentence, or even, an Ehrenfeucht-Fraïssé game of length $\omega$.

It was realised in the 90's by members of the Helsinki school of logic that to develop Scott analysis for uncountable models, one has to play the Ehrenfeucht-Fraïssé game using trees as game clocks (\cite{vaananen1995games}). The trees allow the game to have an arbitrarily long and yet countable length in the same way the ordinal clock $\beta$ makes the game $\EF^\beta$ arbitrarily long but finite.

\begin{de}
    Let $A$ and $B$ be models in a same signature and let $T$ be a tree. The game $\bs{\EF^T(A,B)}$ has two players and is played as follows:
    Positions in the game are pairs $(C,\pi)$, where $C\subseteq T$ and $\pi$ is a function from $A$ to $B$. The starting position is $(\emptyset,\emptyset)$. At position $(C,\pi)$, player $\1$ chooses a node $t'$ and a point $a\in A$ (or a point $b\in B$), and player $\2$ responds by choosing a point $b\in B$ (or a point $a\in A$). The next position is $(C\cup\{t'\},\pi\cup\{(a,b)\})$. 
    The game ends if either the node $t'$ is not $<_T$-above every node in $C$, in which case player $\1$ loses, or else $\pi\cup\{(a,b)\}$ is not a partial isomorphism, in which case player $\2$ loses.

\end{de}

It can first be noted that if the tree $T$ does not have an uncountable branch, then the game $\EF^T$ has countable but arbitrarily long length. Contrary to the game $\EF^\beta$, the game $\EF^T$ might not be determined \cite{hyttinen1987games}. 
Moreover, if $\beta$ is an ordinal and $T_\beta$ is the tree of descending sequences in $\beta$ ordered by end-extension, then the games $\EF^\beta$ and $\EF^{T_\beta}$ are equivalent in the sense that any of the players has a winning strategy in one if and only if they have one also in the other. So the games $\EF^T$ with trees as clocks generalize the games $\EF^\beta$ with ordinals as clocks.

The game with trees as clocks is very effective in measuring the difference between two non-isomorphic models of size $\aleph_1$. To state the analogue of Theorem \ref{thm:scott1}, we first need to introduce an operation, namely \textit{Kurepa's $\sigma$-functor} from \cite{kurepa1956ensembles}. For a tree $T$, the tree $\sigma T$ consists of all the chains in $T$. The tree is ordered by end-extension. Let $S\leq T$ denote the fact that there exists a weak embedding from $S$ to $T$; write $S\equiv T$ if both $S\leq T$ and $T\leq S$ hold and write $S<T$ if $S\leq T$ but $S\not\equiv T$. It can be shown that $T<\sigma T$ and that the $\sigma$-operation extends the successor operation on ordinals: $\sigma T_\alpha\equiv T_{\alpha+1}$ for any ordinal $\alpha$. (In general, it does not hold that if $S<T$, then $\sigma S\leq T$). If $T$ has no uncountable branches, then neither does $\sigma T$. If $T$ is a wide $\aleph_1$-Aronszajn tree and the continuum hypothesis holds, then so is $\sigma T$. However, in general $\sigma T$ might be much wider than $T$. The parameter that measures how different two non-isomorphic models of size $\aleph_1$ are, is a tree without uncountable branches:

\begin{thm}[Hyttinen, Väänänen \cite{hyttinen1990scott}]\label{thm:intro:non-iso} For models $A$ and $B$ of size $\aleph_1$, the following are equivalent:
\begin{enumerate}
    \item $A\not\cong B$,
    \item There are trees $S$ and $T$ without uncountable branches such that $S\leq T$ and
    \begin{enumerate}
        \item player $\2$ has a winning strategy in $\EF^T(A,B)$ but not in $\EF^{\sigma T}(A,B)$,
        \item player $\1$ does not have a winning strategy in $\EF^S(A,B)$ but has one in $\EF^{\sigma S}(A,B)$.
    \end{enumerate}
\end{enumerate}
\end{thm}
The theorem has to be stated in two parts since the Ehrenfeucht-Fraïssé games with trees as games clocks are not necessarily determined \cite{hyttinen1987games}.

While Theorem \ref{thm:intro:non-iso} answers the question how to separate non-isomorphic uncountable models, no complete analogue for Scott's characterization of model up to isomorphism by its Scott sentence (Theorem \ref{thm:scottiso}) has been found. The work is not finished, only partial results exist thus far. In the case of countable models, the logic that corresponds the game $\EF^\beta$ is the logic $\LL_{\infty\omega}$. It is an empirical fact that Ehrenfeucht-Fraïssé types of games always correspond to a logic. In some sense this is the case for the games $\EF^T$ too. There exists a ``tree logic" $\LL_{\TT}$ where trees without uncountable branches work as quantifier ranks. The logic has a generative set of formulas, i.e. the set of formulas is obtained from atomic formulas by closing under a certain set of operations. There is a drawback: the semantics is defined through a variation of a satisfaction game where a tree is used as a clock, which results in a logic that is not necessarily closed under negation, due to the fact that the satisfaction game might not be determined. However, against certain trees the game is always determined, and for some uncountable models $A$ there does exist a tree $T$ such that for any $B$, player $\2$ has a winning strategy in $\EF^T(A,B)$ if and only if $A\cong B$. In this case there also is a Scott sentence $\psi_A$ for $A$ in the tree logic $\LL_{\TT}$. See the survey \cite{vaananen1995games} or \cite{karttunen1984model} and \cite{hyttinen1987games} for more. 

In this paper we show that it is consistent that there exists a universal wide $\aleph_1$-Aronszajn tree (modulo consistency of a weakly compact cardinal). The reason why this scenario is interesting from the point of view of classifying uncountable structures is the following simple observation: if player $\2$ has a winning strategy in the game $\EF^T(A,B)$ and $S\leq T$, then she has a winning strategy also in the game $\EF^S(A,B)$. Thus, in order to verify the similarity of $A$ and $B$ up to a set of potentially incomparable trees $T_i$, $i\in I$, it suffices to verify the similarity of $A$ and $B$ up to a tree that is above any tree $T_i$. Thus, the existence of a universal wide $\aleph_1$-Aronszajn tree acts as a ``master tree" for all the clock trees of size $\aleph_1$.

\vv

\subsection{Preliminaries}\spa

\vv

In this paper we show the consistency of the existence of a universal wide $\aleph_1$-Aronszajn tree. The proof consists of performing a forcing iteration of length $\kappa^+$, where $\kappa$ is a weakly compact cardinal. The first poset collapses the weakly compact onto $\aleph_1$ and introduces a wide $\kappa$-Aronszajn tree $\dot T$, and the rest of the iteration takes care of embedding each wide $\kappa$-Aronszajn tree, also those that appear along the iteration, into $\dot T$. Preservation of $\kappa$ will be guaranteed by strong properness with respect to enough models of size $<\kappa$ and the preservation of the Aronszajn-ness of $\dot T$ will be guaranteed using the reflection properties of the weakly compact. We review basic properties of strongly proper forcing and weakly compact cardinals.

\begin{de}
    Let $\P$ be a poset and let $M$ be a set. 
    \begin{enumerate}
        \item Let $p\in\P$. A condition $r\in\P\cap M$ is a \textbf{residue of $\bs{p}$ into $\bs{M}$} if every $w\in\P\cap M$ that extends $r$ is compatible with $p$.
        \item A condition $p\in\P$ is \textbf{strongly $\bs{(\P,M)}$-generic} if every $q\leq p$ has a residue into $M$.
        \item The poset $\P$ is \textbf{strongly proper with respect to $\bs{M}$} if for every $p\in\P\cap M$ there is $q\leq p$ that is strongly $(\P,M)$-generic.
    \end{enumerate}
\end{de}

\noindent A proof for the following lemma can be found in \cite{gilton2017side}:

\begin{lem}\label{lem:strgen}
    A condition $p$ is strongly $(\P,M)$-generic if and only if
\[
p\Vdash\check G_{\P}\cap M\mbox{ is a }V\mbox{-generic filter on }\P\cap M.
\]
\end{lem} 

\noindent Lemma \ref{lem:strgen} says that if $\P$ is strongly proper with respect to $M$, then $\P\cap M$ is a complete subposet of $\P$ ``modulo a strongly $(\P,M)$-generic condition". Those $V$-generic filters on $\P$ that contain a strongly $(\P,M)$-generic condition project to a generic on $\P\cap M$. Conversely, small generics extend on larger posets, again modulo a strongly generic condition:

\begin{lem}
    Assume that $M$ is a suitable model for $\P$ and $p$ is strongly $(\P,M)$-generic. Then every $V$-generic filter on $\P\cap M$ that contains a residue of $p$ extends to a $V$-generic on $\P$ that contains $p$.
\end{lem}

If $G\subseteq\P\cap M$ is a $V$-generic filter, we define
\[
\P/G:=\{p\in\P:p\text{ is compatible with every }w\in G\}.
\]
Then, a condition $r\in\P\cap M$ is a residue of $p$ into $M$ if and only if $r\Vdash ``p\in\P/\check G_{\P\cap M}"$. It follows from the above lemmas that posets that have strongly generic conditions factor into two, modulo these conditions. Specifically, if $p$ is strongly $(\P,M)$-generic, then any residue $r$ of $p$ forces that $\P/p$ is forcing equivalent to the two-step iteration
\[
(\P\cap M)*(\P/p)/\check G_{\P\cap M}.
\]

It follows from Lemma \ref{lem:strgen} that strong properness is a generalisation of properness. Thus, preservation of cardinals can be expected. The following fact is standard and a proof can be found for instance in \cite{jech2003set}.

\begin{lem}
    Let $\kappa$ be a cardinal. A poset that is strongly proper with respect to stationarily many $M\in\Pow_\kappa(H_\theta)$ for any large enough regular cardinal $\theta$ preserves $\kappa$.
\end{lem}

Let $\kappa\leq\lambda$ be regular cardinals. The Levy collapse poset $\Coll(\kappa,<\lambda)$ is the poset of partial functions $p:\lambda\times\kappa\to\lambda$ of size $<\kappa$ that satisfy $p(\alpha,\beta)<\alpha$, ordered by inverse inclusion.

\begin{ex}\label{ex:coll}
    Let $\kappa\leq \lambda$ be regular cardinals and let $M$ be a model of $\ZFC^-$ of size $<\lambda$ that satisfies $\lambda\cap M\in\Ord$ and contains $\kappa$ as element. Then the Levy Collapse $\Coll(\kappa,<\lambda)$ is strongly proper with respect to $M$. Then the intersection $p\cap M$ is a residue of $p$ into $M$ for any condition $p\in\Coll(\kappa,<\lambda)$.
\end{ex}

\noindent For more about the general theory of strong properness see for example \cite{gilton2017side}.

\vv

A formula $\phi$ is a \textbf{$\bs{\Pi^1_1}$-formula over a set $\bs{M}$} if $\phi$ has form $\forall X\psi(X,\bar{A},\bar{a})$, where $\bar{A}=(A_1,\dots,A_n)$ is a tuple of subsets of $M$ and $\bar{a}=(a_1,\dots a_m)$ is a tuple of elements of $M$, and $X$ is a second-order variable. Second-order variables are interpreted as unary predicates.

\begin{de}
    An uncountable cardinal $\kappa$ is \textbf{weakly compact} if for every $\Pi^1_1$-formula $\phi$ and $A\subseteq V_\kappa$ such that 
    \[
    V_\kappa\models\phi(A),
    \]
    there is $\alpha<\kappa$ such that
    \[
    V_\alpha\models\phi(A\cap V_\alpha).
    \]
\end{de}

If $\kappa$ is weakly compact, then the sets of the form
\[
X_{A,\phi}=\{\alpha<\kappa:\text{if } V_\kappa\models\phi\text{ then }V_\alpha\models\phi(A\cap V_\alpha)\}
\]
generate a filter on $\kappa$. This filter is called the \textbf{weakly compact filter} and denoted by $\bs{\FF_{\mathsf{wc}}}$. It is normal and $\kappa$-complete, and it extends the club filter. See chapter 1 section 6 from \cite{kanamori2008higher} for proofs.

The following observations are used tacitly in the last section of this paper. They can be proved by applying Lemma 13.12 from \cite{jech2003set} to the recursive definition of the forcing relation.

\begin{obs}
    Let $\kappa$ be a cardinal. If $\phi$ is a $\Pi^1_1$-formula over $V_\kappa$ and $\P\subseteq V_\kappa$ is a poset, then so is the formula
    \[
    p\Vdash_{\P}\phi.
    \]
    Furthermore, for every $\P$-name there is an equivalent $\P$-name that is a subset of $V_\kappa$.
\end{obs}

\begin{obs}
    If $\P$ is a poset such that $\P\subseteq V_\kappa$, then the formula
    \[
    \spa ``p\Vdash_\P\dot T\text{ is a wide }\kappa\text{-Aronszajn tree}"
    \]
    is a $\Pi^1_1$-formula over $V_\kappa$.
\end{obs}

\section{The poset}

\noindent This section is devoted to defining the forcing poset that will create the universal wide $\aleph_1$-Aronszajn tree. Throughout, we assume that $\GCH$ holds in $V$.

\subsection{Adding a flexible tree}\spa 

\vv

Let $\kappa$ be a regular uncountable cardinal. In this subsection we define a $\Coll(\omega,<\kappa)$-name $\dot T$ for a wide $\kappa$-Aronszajn tree, which will later, under further forcing, be made universal for wide $\kappa$-Aronszajn trees. Later the cardinal $\kappa$ will be assumed to be weakly compact.

Recall that a wide $\kappa$-tree is a tree of size $\kappa$ and a tree is normal if the meet of two nodes is always well-defined and unique.

\begin{prop}\label{prop:1tree}
    Let $\kappa$ be an uncountable regular cardinal. There is a $\Coll(\omega,<\kappa)$-name $\dot T$ for a tree such that the following are satisfied:
    \begin{enumerate}
        \item In $V^{\Coll(\omega,<\kappa)}$, $\dot T$ is a normal wide $\kappa$-Aronszajn tree such that any node $t\in T$ has $\kappa$ many successors at any level above. 
        \item The domain of $\dot T$ is the set $T:=\kappa\times\kappa$ and $\Lev_\alpha(\dot T)=\kappa\times\{\alpha\}$ for each $\alpha<\kappa$.
        \item There are club many $\alpha<\kappa$ such that for every generic filter $G\subseteq\Coll(\omega,<\kappa)$,
        \[
        \dot T^G\cap V_\alpha\in V[G\cap V_\alpha].
        \]
        \item\label{prop:1tree:firstitems} There is a function $\width:T\to\kappa$ in $V$ such that for any node $t\in T$, $\width(t)$ is forced by $\emptyset$ to be the least ordinal $\width(t)\geq\height(t)$ such that the branch below $t$ is contained in $\width(t)\times\height(t)$. For any $\alpha\leq\beta$, there are $\kappa$ many nodes $t\in\Lev_\alpha(T)$ with $\width(t)=\beta$.
        \item\label{prop:1tree:pi} \emph{(Projections $\pi^p(t)$.)} For every node $t\in T$ and a condition $p\in\Coll(\omega,<\kappa)$, the set 
        \[
        \{t'\in T:p\Vdash ``t'<_{\dot T}t"\}
        \]
        has an element at a maximum height. We denote it by $\bs{\pi^p(t)}$ and call it the \textbf{node projection of $\bs{t}$ with respect to $\bs{p}$}. 
        \item\label{prop:1tree:free} \emph{(Flexibility.)} For $t\in T$, we say that a condition $p\in\Coll(\omega,<\kappa)$ is \textbf{nice with respect to $\bs{t}$}, if for any $\alpha<\kappa$ with $t\in T-V_\alpha$, any condition $w\in\Coll(\omega,<\kappa)\cap V_\alpha$ and any node $\bar{t}\in\width(t)\times\height(t)$, if $w\leq p\cap V_\alpha$ and
        \[
        w\Vdash``\pi^p(t)<_{\dot T}\bar{t}",
        \]
        then there is a condition $q\leq p$ such that
        \begin{enumerate}
            \item $q\cap V_\alpha=p\cap V_\alpha$,
            \item any common extension of $w$ and $q$ forces $``\bar{t}<_{\dot T}t"$,
            \item $q$ is a minimal extension of $p$ in the sense that if $t'\neq t$, then
            \[
            \pi^q(t')=\pi^p(t').
            \]
        \end{enumerate}
        We require that for any set of nodes $A\subseteq T$, the conditions that are nice with respect to every $t\in A$ are dense in $\Coll(\omega,<\kappa)$.
        \item\label{item:labels:1tree} \emph{(Labels.)} There is a function $\lab:T\to\kappa$ with the following properties:
        \begin{enumerate}
            \item $\width(t)\leq\lab(t)$ for every $t\in T$,
            \item for any $\alpha\leq\beta_0\leq\beta_1$, there are $\kappa$ many nodes $t\in\Lev_\alpha(T)$ with $\width(t)=\beta_0$ and $\lab(t)=\beta_1$.
        \end{enumerate}
        The ordinal $\bs{\lab(t)}$ is called the \hypertarget{label}{\textbf{label of $\bs{t}$}}. 
\end{enumerate}
\end{prop}
\begin{proof} We build a $\Coll(\omega,<\kappa)$-name $\dot T$ for a wide $\kappa$-Aronszajn tree on a set of size $\kappa$. The domain $T$ of the tree $\dot T$ will be the set $\kappa\times\kappa$ and each $\alpha$th level will be the set $\kappa\times\{\alpha\}$. 
The tree-order $<_{\dot T}$ on $T$ will be defined in $V^{\Coll(\omega,<\kappa)}$. First, we fix two functions $\width:T\to \kappa$ and $\lab:T\to\kappa$ (in $V$) with the following properties:
\begin{enumerate}
    \item for any node $t=(\beta_t,\alpha_t)\in T$, if $\beta_t\leq\alpha_t$, then $\height(t)=\alpha_t=\width(t)=\lab(t)$, and if $\alpha_t<\beta_t$, then 
    \[
    \height(t)=\alpha_t\leq\width(t)\leq\lab(t)\leq\max\{\alpha_t,\beta_t\},
    \]
    \item for any $\alpha\leq\beta_0\leq\beta_1<\kappa$, there are $\kappa$ many nodes $t\in\Lev_\alpha(T)$ with $\width(t)=\beta_0$ and $\lab(t)=\beta_1$. 
\end{enumerate}
Then, denote by $\Q$ the poset consisting of finite partial functions $p:T\times \omega\to T$ such that if $(t,n)\in\dom(p)$, then
\[
f^p_t(n):=p(t,n)\in\width(t)\times\height(t),
\]
ordered by inverse inclusion. A generic filter $G\subseteq\Q$ yields, for every $t\in T$, a surjection
\[
f^G_t:\omega\to\width(t)\times\height(t).
\]
It can be seen that $\Q$ is isomorphic to $\Coll(\omega,<\kappa)$. We are ready to define the tree order $<_{\dot T}$. By the isomorphism, it suffices to define a $\Q$-name for it, instead of $\Coll(\omega,<\kappa)$-name. To this end, let $G\subseteq\Q$ be a generic filter.
Working in $V[G]$, we define the order $<_{T}$ on $T$ by recursion on levels $\alpha<\kappa$. Suppose that the order has already been defined on $\Lev_{<\alpha}(T)$.

For every $t\in \Lev_\alpha(T)$ we define the branch below it, denoted by $b_t$. This suffices, since then we let $s<_Tt$ if $s$ belogns to the branch $b_t$. The branch $b_t$ is derived from the function
\[
f^G_t:=\bigcup_{p\in G}p(t,\cdot):\omega\to \width(t)\times\height(t),
\]
by recursion on $m<\omega$ as follows: let $t_0$ be $f^G_t(0)$, and whenever $t_m$ is defined, look at the least $n\in[m,\omega)$ such that the node $f^G_t(n)$ is strictly above $t_m$ in the tree order $(\Lev_{<\alpha}(T),<_T)$ and let $t_{m+1}:=f^G_t(n)$.
The branch below $t$ can now be defined to be the set
\[
b_t:=\bigcup_{m<\omega}\{s\in\Lev_{<\alpha}(T): s<_T t_m\}.
\]
It follows by a density argument that $b_t$ is a cofinal branch through $\Lev_{<\alpha}(T)$. This concludes the definition of the tree $(T,<_T)$. Let $<_{\dot T}$ be a $\Q$-name for the order $<_T$ and denote $\dot T=(T,<_{\dot T})$. It remains to verify that $\dot T$ satisfies the proposition.

\begin{claim}\label{prop:1tree:claim:st}
    The tree $\dot T$ satisfies the properties of the proposition.
\end{claim}
\begin{proof}[Proof of Claim \ref{prop:1tree:claim:st}] \spa

\textbf{Items (1)-(\ref{prop:1tree:firstitems}).} It is clear that in $V^{\Q}$,  $\dot T$ is a wide $\kappa$-tree. For distinct nodes $t$ and $t'$, the branches $b_t$ and $b_{t'}$ are eventually disjoint, since the functions $f^G_t$ and $f^G_{t'}$ must be mutually generic. Thus $\dot T$ is normal. The choice of the width function $\width$ implies that every node has $\kappa$ many successors at any level above. Whenever $V_\alpha$ is closed for both the label function $\lab$ and the width function $\width$, then $\dot T^G\cap V_\alpha$ is a subtree of $\dot T^G$ and $\dot T^G\cap V_\alpha\in V[G\cap V_\alpha]$. There are club many such $\alpha$'s. Furthermore, by a density argument, the ordinal $\width(t)$ must be the least such that the branch below $t$ is contained in the rectangle $\width(t)\times\height(t)$. We postpone the proof that $\dot T$ has no cofinal branches.

\textbf{Item (\ref{prop:1tree:pi}).} Let $t\in T$ and $p\in\Q$. The set 
\[
\{t'\in T:p\Vdash t'<_{\dot T}t\}
\]
is contained in the set $\ran(p_t)$, which is finite. Hence it has a node at a maximum level. This set is linearly ordered, hence it has a unique maximum.

\textbf{Item (\ref{prop:1tree:free}).} Let $p$ be a condition. We say that $p$ satisfies $\bs{*(t)}$ if the domain of the function $f^p_t$ is a natural number and for any $t_0,t_1\in\ran(f^p_t)$, $p$ decides whether $t_0<_{\dot T}t_1$ or $t_0\not<_{\dot T}t_1$. For any set $t\in T$, we may extend the functions $f^p_{t_1}$ for $t_1\in\ran(f^p_t)$. Only finitely many functions must be extended to make $p$ satisfy $*(t)$. Similarly, it can be seen that for any $A\subseteq T$, the set of conditions that satisfy $*(t)$ for every $t\in A$ is dense in $\Q$. 

We show that if $p$ satisfies $*(t)$, then it is nice with respect to $t$. Let $t\in T$ let $p$ satisfy $*(t)$. Let $\alpha<\kappa$ be such that $t\in T-V_\alpha$, and let $w\in\Q\cap V_\alpha$ and $\bar{t}\in\width(t)\times\width(t)$ be such that $w\leq p\cap V_\alpha$ and
\[
w\Vdash``\pi^p(t)<_{\dot T}\bar{t}".
\]
Suppose that $\dom(f^p_t)=n$. Then there is $m<n$ such that $p\Vdash ``\pi^p(t)=t_m"$, where $t_m$ is as in the definition of the branch below $t$, and for every $l\in [m,n)$, we have $p\Vdash``\pi^p(t)\not<_{\dot T}f^p_t(l)$, by $*(t)$ and maximality of $\pi^p(t)$. Define $q$ from $p$ by letting $f^q_t:=f^p_t\cup\{(n,\bar{t})\}$ and $f^q_{t'}:=f^p_{t'}$ for $t'\neq t$. Then $q\cap V_\alpha=p\cap V_\alpha$ because $t\in T-V_\alpha$. It is also clear that $\pi^q(t')=\pi^p(t')$ for any $t'\neq t$. We claim that any common extension of $w$ and $q$ forces $\bar{t}<_{\dot T}t$. To see this, fix a generic $G\subseteq\Q$ that contains the condition $q\cup w$. Then $t_m=\pi^p(t)$ and $t_{m+1}=\bar{t}$, since $n$ is the least above $m$ such that $t_m<_{\dot T^G}f^q_t(n)$, and $f^q_t(n)=\bar{t}$. This is enough to show that nodes that are nice with respect to nodes in $A$ for any $A\subseteq T$ are dense in $\Q$.

\textbf{$\boldsymbol{\dot T}$ has no cofinal branches.} We show that $\dot T$ is a name for a wide $\kappa$-Aronszajn tree. Suppose to the contrary that there is a condition $p\in\Q$ and a name $\dot b$ such that
\[
p\Vdash``\dot b\text{ is a cofinal branch through }\dot T".
\]
Let $\theta$ be some large enough cardinal, let $M\elem H_\theta$ be a model which is closed under the function $t\mapsto R(t)$, such that $\alpha:=\kappa\cap M\in\kappa$ and $p,\dot b\in M$. Find $q\leq p$ that decides the node $t:=\dot b(\alpha)$ and belongs to the dense set from item (\ref{prop:1tree:free}). Then find $w\leq q\rest\alpha$ in $\Q\cap M$ and two distinct nodes $t^L$ and $t^R$ in $\width(t)\times\height(t)$ at the some height $\bar{\alpha}<\alpha$ such that
\[
w\Vdash ``\pi^q(t)<_{\dot T}t^L,t^R".
\]
By item (\ref{prop:1tree:free}) there are two extensions $q^L,q^R\leq q$ such that $q^L\rest\alpha=q^R\rest\alpha=q\rest\alpha$ and such that any common extension of $w$ and $q^L$ (respectively $q^R$) forces $``t^L<_{\dot T}t"$ (respectively $``t^R<_{\dot T}t"$).
Up to extending $w$ further inside $M$, we may assume that it decides the node $\bar{t}:=\dot b(\bar{\alpha})$. If $\bar{t}\neq t^L$, then $w\cup q^L$ is a condition that forces both $\bar{t}<_{\dot T}t$ and $t^L<_{\dot T}t$, which is contradictory, and if $\bar{t}\neq t^R$, then $r\cup q^R$ is a condition that forces both $\bar{t}<_{\dot T}t$ and $t^R<_{\dot T}t$, which is also contradictory. Hence $\dot T$ must be a wide $\kappa$-Aronszajn tree.

\end{proof}

\noindent This ends the proof of the proposition.
    
\end{proof}

\subsection{Adding embeddings}\spa

\begin{assumption*}
Henceforth, we fix a weakly compact cardinal $\kappa$, a $\Coll(\omega,<\kappa)$-name for a wide $\kappa$-Aronszajn tree $\dot T$ and width and labeling functions $t\mapsto\width(t)$ and $t\mapsto \lab(t)$, as in Proposition \ref{prop:1tree}. 
\end{assumption*}

\begin{notation}
    We say that a collapse condition $p\in\Coll(\omega,<\kappa)$ is \textbf{nice} if it belongs to the dense set of item (\ref{prop:1tree:free}) of Proposition \ref{prop:1tree}.
\end{notation}


We will define a finite support iteration $(\P_\delta:\delta\leq\kappa^+)$ such that the first poset is the poset $\Coll(\omega,<\kappa)$ and the final poset $\P_{\kappa^+}$ adds an embedding from every wide $\kappa$-Aronszajn tree into $\dot T$.
We begin by defining sequences of suitable models of size $<\kappa$, which will be used as side conditions. We will define by recursion on $\delta<\kappa^+$, in this order: the poset $\P_\delta$, models $(M^\delta_\alpha:\alpha<\kappa)$ and a set $\EE_\delta\subseteq\kappa$ that selects the ``good models" $(M^\delta_\alpha:\alpha\in\EE_\delta)$. These models will then be used in the definition of the poset $\P_{\delta+1}$.

\begin{notation}
    Fix a large enough regular cardinal $\theta$ and a well-order $<_\theta$ of $H_\theta$.
\end{notation}

The \textbf{weakly compact filter $\bs{\FF_{\mathsf{wc}}}$} is the normal filter on $\kappa$ generated by sets $\{\alpha<\kappa:$ if $V_\kappa\models\phi(A)$, then $V_\alpha\models\phi(A\cap V_\alpha)\}$, for every $\Pi^1_1$-formula $\phi$ and set $A\subseteq V_\kappa$.

\begin{de}
    Let $\alpha<\kappa$. A model $M\elem H_\theta$ \textbf{\hypertarget{reflects}{reflects} at $\boldsymbol{\alpha}$} if for every $A\in\Pow(V_\kappa)\cap M$ and $\Pi^1_1$-formula $\phi$,
    \[
    V_\kappa\models\phi(A)\implies V_\alpha\models\phi(A\cap V_\alpha).
    \]
\end{de}

\begin{de}\label{de:Edelta}
    Let $\P_0:=\{\emptyset\}$. By recursion on $\delta<\kappa^+$: whenever the poset $\P_\delta$ and the sets $\bar{\EE}_\delta:=(\EE_\gamma)_{\gamma<\delta}$ are defined:
    \begin{enumerate}
        \item for every $\alpha\leq\kappa$, let \[M^\delta_\alpha:=\text{Skolem hull of }\alpha\cup\{\kappa,\dot T,\delta,\P_\delta,\bar{\EE}_\delta\}\text{ in }(H_\theta,\in,<_\theta),\]
    \item let
        \begin{align*}
            \hypertarget{EE}{\EE_\delta}:=\{\alpha<\kappa:\spa (1)\quad & V_\kappa\cap M^\delta_\alpha=V_\alpha,\\
            (3)\quad &M^\delta_\alpha\text{ \hyperlink{reflects}{reflects} at }\alpha, \\
            (4)\quad &\alpha\in\bigcap_{\gamma\in\delta\cap M^\delta_\alpha}\lim\EE_\gamma\}.
        \end{align*}
    \end{enumerate}
\end{de}

As mentioned above, the models $(M^\delta_\alpha:\alpha\in\EE_\delta)$ are the ones that will be used as side conditions when defining the poset $\P_{\delta+1}$. We assume that each $M^\delta_\alpha$ is closed for the width and labeling maps $t\mapsto \width(t)$ and $t\mapsto \lab(t)$.

\begin{lem}\label{lem:structureofmodels} \spa 

\begin{enumerate}
    \item The set $\EE_\delta$ belongs to the weakly compact filter $\FF_{\mathsf{wc}}$ for every $\delta<\kappa^+$.
    \item For every $\delta<\kappa^+$ and $\alpha\in\EE_\delta$:
    \begin{enumerate}
        \item $\P_\delta\in M^\delta_\alpha$,
        \item $\P_\gamma\in M^\delta_\alpha$ for every $\gamma\in\delta\cap M^\delta_\alpha$,
        \item\label{lem:structureofmodels:item3} $M^\gamma_\beta\in M^\delta_\alpha$ for every $\beta<\alpha$ and $\gamma\in M^\delta_\alpha$,
        \item\label{lem:strmod:item5} $\delta\cap M^\delta_\alpha=\delta\cap M^{\delta+1}_\alpha$ for all $\alpha<\kappa$ and $\delta<\kappa^+$.
        \item\label{lem:structureofmodels:item4} if $\xi<\gamma<\delta$ and $\xi,\gamma\in M^\delta_\alpha$, then $\xi\in M^\gamma_\alpha$,
    \end{enumerate}
\end{enumerate}
\end{lem}
\begin{proof}
\spa 

\begin{enumerate}
    \item This follows from the assumption that $\kappa$ is weakly compact and the normality of the filter $\FF_{\mathsf{wc}}$.
    \item For item (\ref{lem:structureofmodels:item3}): if $\beta<\alpha$ and $\gamma\in M^\delta_\alpha$, then $M^\gamma_\beta\subseteq M^\delta_\alpha$, and since $M^\delta_\alpha$ is closed under sequences of length $<\alpha$ and $|M^\gamma_\beta|=\beta<\alpha$, we obtain $M^\gamma_\beta\in M^\delta_\alpha$.
    For item (\ref{lem:strmod:item5}): if $\delta<\kappa$, then $\delta<\alpha$ and the claim follows easily, and otherwise if $\psi:\kappa\to\delta$ is the $<_\theta$-least bijection, then $\delta\cap M^{\delta+1}_\alpha=\psi[\alpha]=\delta\cap M^\delta_\alpha$.
    Item (\ref{lem:structureofmodels:item4}) is proved similarly: if $\gamma<\kappa$, then $\gamma<\alpha$ and the claim follows easily, and otherwise if $\psi:\kappa\to\gamma$ is the $<_\theta$-least bijection, then $\psi\in M^\delta_\alpha$ and $\psi\in M^\gamma_\alpha$, which implies $\gamma\cap M^\delta_\alpha=\psi[\alpha]=\gamma\cap M^\gamma_\alpha$, and gives directly $\xi\in M^\gamma_\alpha$.

    \end{enumerate}
    
\end{proof}

We next fix a bookkeeping function. We intentionally leave some flexibility regarding it, for instance, for now we will only assume that it picks trees of size $\kappa$ which are not necessarily assumed to be Aronszajn.

\begin{notation}\label{nota:bookkeeping} 
    Fix a bookkeeping function \[
    \kappa^+\to \Pow(V_\kappa),\quad\gamma\mapsto\dot S_\gamma
    \] such that whenever the poset $\P_\gamma$ has been defined, then $\dot S_\gamma$ is a $\P_\gamma$-name for a normal tree with domain $\kappa\times\kappa$ and that $\Vdash_{\P_\gamma}\Lev_\alpha(\dot S_\gamma)=\kappa\times\{\alpha\}$. We assume without loss of generality that each model $M^\delta_\alpha$ contains the bookkeeping function. We denote by $S_\gamma$ the domain of $\dot S_\gamma$, i.e. the set $\kappa\times\kappa$.
\end{notation}

\begin{de}
    Let $S$ be a tree. A node $s\in S$ is an \textbf{exit node} from a set $M$ if $s\notin M$ but $b_s\subseteq M$, i.e. the branch below $s$ is contained in $M$.
\end{de}

We will now define the posets $\P_\delta$, $\delta\leq\kappa^+$. Inductively, we assume that the sets \hyperlink{EE}{$\EE_\gamma$, ${\gamma<\delta}$} from Definition \ref{de:Edelta}, have already been defined, when defining $\P_\delta$.

\begin{de}\label{de:poset} Let $\delta\leq\kappa^+$. Conditions in the poset $\P_\delta$ are functions
\[
p:\delta\to V_\kappa
\]
such that for every $\gamma<\delta$,
\[
p(\gamma)=(f^p_\gamma,\NN^p_\gamma)
\]
satisfies
\begin{enumerate}
    \item $f^p_0\in\Coll(\omega,<\kappa)$,
    \item for non-zero $\gamma<\delta$, $f^p_\gamma:\kappa\times\kappa\to \kappa\times\kappa$ is a finite partial injective map such that:
    \begin{enumerate}
        \item $p\rest\gamma$ decides the $\dot S_\gamma$-meets in the set $\dom(f^p_\gamma)$,
        \item $p\rest\gamma\VVdash{\dom(f^p_\gamma)\mbox{ is closed under }\dot S_\gamma\text{-meets}}$,
        \item $p\rest\gamma\VVdash{f^p_\gamma:\dot S_\gamma\to \dot T\text{ is a level- and meet-preserving tree-embedding}}$,
    \end{enumerate}
    \item $\NN^p_\gamma\subseteq\lim\EE_\gamma$ is a finite set such that whenever $\alpha\in \NN^p_\gamma$ and $\xi\in\gamma\cap M^\gamma_\alpha$, then $\alpha\in \NN^p_\xi$, and such that the union
    \[
    \bigcup_{\gamma\in\delta}\NN^p_\gamma
    \]
    is finite,
    \item for every non-zero $\gamma<\delta$, every $s\in\dom(f^p_\gamma)$ and $\alpha\in\NN^p_\gamma$:
    \begin{enumerate}
        \item $s\in M^\gamma_\alpha$ if and only if $f^p_\gamma(s)\in M^\gamma_\alpha$,
        \item $p\rest\gamma\VVdash{s\mbox{ is an exit node from }M^\gamma_\alpha\mbox{ if and only if }f^p_\gamma(s)\mbox{ is}}$,
        \item if $p\rest\gamma\VVdash{s\mbox{ is an exit node from }M^\gamma_\alpha}$, then $\width(f^p_\delta(s))\in\NN^p_\gamma$ and $p\rest\gamma$ decides the ordinal 
        \[
        \width(s):=\text{the least }\beta\geq\height(s)\text{ such that }b_s\subseteq\beta\times\height(s),
        \] 
        and it satisfies $\width(s)\leq \width(f^p_\delta(s))$, 
        \item if $s\notin M^\gamma_\alpha$, then there is $\bar{s}\in\dom(f^p_\gamma)-M^\gamma_\alpha$ such that
        \[
        p\rest\gamma\VVdash{\bar{s}<_{\dot S_\gamma}s\text{ and }\bar{s}\mbox{ is an exit node from }M^\gamma_\alpha},
        \]
        \item\label{de:posetseparation} if $p\rest\gamma\VVdash{s\mbox{ is an exit node from }M^\gamma_\alpha}$, then the label $\lab(f^p_\gamma(s))$ satisfies:
        \begin{enumerate}
            \item $\lab(f^p_\gamma(s))\in \EE_\gamma$,
            \item $s\in M^\gamma_{\lab(f^p_\gamma(s))}$, 
            \item $\lab(f^p_\gamma(s))\in \NN^p_\xi$ for every $\xi\in\gamma\cap M^\gamma_{\lab(f^p_\gamma(s))}$,
        \end{enumerate}
        \end{enumerate}
        \item\label{de:posetshrink} the support $\sp(p):=\{\gamma\in\delta:f^p_\gamma\neq\emptyset\}$ is finite.
    \end{enumerate}
The order is the pointwise inverse inclusion: $q\leq p$ if $f^q_\gamma\supseteq f^p_\gamma\mbox{ and }\NN^q_\gamma\supseteq\NN^p_\gamma$ for all $\gamma<\delta$.
\end{de}

\begin{rmk}
    If $\rho$ is the label of $t$ and $\rho\in\EE_\delta$, then $t\notin M^\delta_\rho$ and $t$ is an exit node from $M^\delta_\rho$.
\end{rmk}

\begin{lem} If $\gamma<\delta$ and $q\in\P_\gamma$ extends $p\rest\gamma$, then the concatenation $q^\smallfrown p\rest[\gamma,\delta)$ is a condition in $\P_\delta$ extending $q$ and $p$. In particular $\P_\gamma\subseteq_c\P_\delta$.    
\end{lem}

\begin{notation}
    Each $\P_\delta$, $\delta<\kappa^+$, has size $\kappa$ and can thus be coded as a subset of $V_\kappa$ using any injection $\delta\to\kappa$. We tacitly identify each $\P_\delta$ with an isomorphic poset contained in $V_\kappa$, and assume that this identification was done using a coding function belonging to each model $M^\delta_\alpha$, $\alpha<\kappa$.
\end{notation}

\begin{rmk}
    The set $\{\gamma<\delta:\NN^p_\gamma\neq\emptyset\}$ is not necessarily finite, even if the union $\bigcup_{\gamma<\delta}\NN^p_\gamma$ is.
\end{rmk}

\begin{rmk}\label{rmk:cc}
    The final poset $\P_{\kappa^+}$ is the direct limit of the posets $(\P_\delta:\delta<\kappa^+)$, via the maps $p\mapsto p^\smallfrown(\emptyset,\emptyset),\dots$. In particular, it follows that $\P_{\kappa^+}$ has $\kappa^+$-cc, being a direct limit of posets of size $\kappa$.
\end{rmk}

The iteration $(\P_\delta:\delta\leq\kappa^+)$ has now been defined. It will be seen that it has the following properties:
\begin{enumerate}
    \item $\P_{\kappa^+}$ has $\kappa^+$-cc, so it preserves all cardinals $\lambda\geq\kappa^+$,
    \item $\P_{\kappa^+}$ collapses every $\alpha<\kappa$ to be countable,
    \item\label{item:presofkappa} $\P_{\kappa^+}$ preserves $\kappa$ and adds $\kappa^+$ many reals,
    \item there is an injective level-preserving tree-embedding $\dot f_\gamma:\dot S_\gamma\to\dot T$ in $V^{\P_{\kappa^+}}$, for every $\gamma<\kappa^+$,
    \item if the bookkeeping function picks only names for wide $\kappa$-Aronszajn trees, then $\dot T$ is a wide $\aleph_1$-Aronszajn tree in $V^{\P_{\kappa^+}}$.
\end{enumerate}
In particular, with a suitable bookkeeping function, $\dot T$ will be a universal wide $\aleph_1$-Aronszajn tree in $V^{\P_{\kappa^+}}$.

The first two properties are immediate from Remark \ref{rmk:cc} and the definition of the poset.
The rest of the paper is devoted to proving the remaining three properties.

\section{Properties of the poset}

In this section, it will be shown that $\P_\delta$, for each $\delta<\kappa$, is strongly proper with respect to stationarily many $M\in[H_\lambda]^{<\kappa}$, for any large enough regular cardinal $\lambda$, and furthermore, that for every $\gamma<\kappa^+$ and $s\in S_\gamma$, the set of conditions $p$ with $s\in\dom(f^p_\gamma)$ is dense. This will imply that $\P_{\kappa^+}$ preserves $\kappa$, making it the new $\aleph_1$, that $\P_{\kappa^+}$ adds $\kappa^+$ many reals, making continuum to be of size $\aleph_2$, and that $\dot T$ is made universal for wide $\kappa$-Aronszajn trees. The section after this one is dedicated to showing that $\dot T$ itself remains (wide) Aronszajn.

\vv 

We now decribe informally the simple idea behind the proof of strong properness. First, note why it is important that the tree $\dot T$ is not in $V$, but rather introduced by the first poset $\P_1$. Indeed, suppose for the sake of argument that we tried to add a tree-embedding 
\[
\dot f: S\to T
\]
using finite approximations $p:S\to T$ of a tree-embedding as conditions, from a (wide) $\aleph_1$-tree $S$ into a tree $T$ both of which are in $V$. The ordering is inverse inclusion. Suppose we now want to show that this poset is strongly proper. Fix some countable model $M$ such that $S,T\in M$. We argue that no condition $p$ can have a residue into $M$. Indeed, fix some condition $p$ which is strong enough that it contains some nodes in $S$ that are in $S-M$ and suppose that it had a residue $r$ into $M$. Now choose a node $s\in S-M$ at height $\omega_1\cap M$ whose predecessors are all in $M$. Look at the branch below $s$ - there is a maximal node $\bar{s}$ which is both below $s$ and in the domain of $r$. Choose some node $s'$ which is above $\bar{s}$ and below $s$ and another node $t'$ in $T\cap M$ which is above $r(\bar{s})$ but \textit{not} below $p(s)$. Extend the condition $r$ by letting $w:= r\cup\{(s',t')\}$. Then $w$ extends $r$ and belongs to $M$. However, this $w$ cannot be compatible with the original condition $p$ since a tree-embedding must be order-preserving. Any common extension of $w$ and $p$ would map the node $s'$ which is below $s$ to the node $t'$ which is not below $p(s)$.

We solve this problem by taking the tree $T$ to be generic over $V$. Indeed, suppose now that our conditions are pairs $p=(u^p,f^p)$, where $u^p$ is a finite approximation of a wide $\aleph_1$-tree $\dot T$ and $u^p$ forces that $f^p$ is a finite approximation of a tree-embedding
\[
f^p:S\to\dot T,
\]
where $S$ is some $\aleph_1$-tree in $V$. We show that the above problem does not occur. Suppose for simplicity that we are in the following situation: $p$ is a condition and $M$ is a model such that the model $M$ is closed under the finite tree-embedding $f^p$, i.e. a node $s\in\dom(f^p)$ is in $M$ if and only if its image $f^p(s)$ is in $M$ (in our case, this situation will be guaranteed by the usage of side conditions). Now consider the trace $p\cap M:=(u^p\cap M,f^p\cap M)$. This is still a condition. For every node $t\in \ran(f^p)-M$, there is a node $\pi(t)$ which is forced by $u^p$ to be the maximal node which is both below $t$ and in $M$. Namely, this $\pi(t)$ is the node in $u^p\cap M$ which is below $t$ in the order $<_{u^p}$ and at a maximal height. Now, we extend the trace $p\cap M$ adding a preimage for the node $\pi(t)$. Whenever $(s,t)\in f^p-M$, we look at the predecessor of $s$ at the height if $\pi(t)$ and call it $\bar{s}$. We then add the pairs $(\bar{s},\pi(t))$ to the trace. In other words, we define
\[
r:=(p\cap M)\cup\{(\bar{s},\pi(t)):(s,t)\in f^p-M\}.
\]
It can be shown that $r$ is a residue of $p$ into $M$. For instance, the previous problem cannot arise: suppose for simplicity that $w$ extends $r$ just by one pair $(s',t')$ and suppose that $s'$ is below some node $s\in\dom(f^p)-M$. But now the condition $p$ does not decide whether $t'$ is below $f^p(s)$ or not. Therefore, we are free to take the pointwise union of $w$ and $p$ and extend its tree-part $u^w\cup u^p$ to some $u$ which decides that $t'$ is indeed below $f^p(s)$. If it happened that $s'$ was not below $s$, we would have done the same, except extended the tree part $u^w\cup u^p$ to some $u$ that decides that $t'$ is not below $f^p(s)$. This is the idea in the proof of strong properness. The difference is that we need to take care of trees $\dot S$ that are not necessarily in $V$, and we have multiple embeddings to deal with. This idea will be exploited in what follows.

\vv 

We begin with a version of a density lemma that will be used in the proof of strong properness. It is conditioned on the strong properness of $\P_\delta$. The lemma will be used in the proof of strong properness of $\P_{\delta+1}$.

\vv

\noindent We need the following notation.

\begin{de}[$\EE^p_\delta$]
    Let $\delta<\kappa^+$ and $p\in\P_\delta$. Define
    \[
    \EE^p_\delta:=\{\alpha\in\EE_\delta:\alpha\mbox{ belongs to }\NN^p_\gamma\mbox{ for every }\gamma\in\delta\cap M^\delta_\alpha\}
    \]
\end{de}

\subsection{Node density}\spa

\begin{lem}[Node Density]\label{lem:dens} Let $\delta<\kappa^+$. Assume that for every condition $p\in\P_\delta$, if $\alpha\in\EE^p_\delta$, then $p$ has a residue into $M^\delta_\alpha$. For every $p\in\P_{\delta+1}$ and $s\in S_\delta$ there is $q\leq p$ such that $s\in\dom(f^q_\gamma)$ and moreover, the extension $q\leq p$ is minimal in the sense that $\NN^q_\delta=\NN^p_\delta$ and $q$ forces that $\dom(f^q_\delta)$ is the least set containing $\dom(f^p_\delta)\cup\{s\}$ and being closed under $\dot S_\delta$-meets and exit nodes from $M^\delta_\alpha$ for $\alpha\in\NN^q_\delta$.
\end{lem}

\begin{proof} 
Let $p\in\P_{\delta+1}$ and $s\in S_\delta$. Up to extending $p\rest\delta$, we may assume that it decides the non-trivial meets in $\dom(f^p_\delta)\cup\{s\}$. In fact, there is at most one meet $s\wedge s'$ not already in the domain of $f^p_\delta$. We may also assume that $p\rest\delta$ decides its implicit image
\[
t_{s'\wedge s}:=\text{the unique node below }f^p_\delta(s')\text{ at the height of }s'\wedge s,
\]
where $s'\in\dom(f^p_\delta)$ is any node such that $s\wedge s'\notin\dom(f^p_\delta)$.
It is starightforward to verify that \[
{p\rest\delta}^\smallfrown(f^p_\delta\cup\{(s'\wedge s,t_{s'\wedge s})\},\NN^p_\delta)
\] is a condition in $\P_{\delta+1}$ extending $p$. Therefore from now onwards, we will assume without loss of generality that $s'\wedge s\in\dom(f^p_\delta)$ and only worry about $s$ itself and exit nodes below $s$ from models $M^\delta_\alpha$, where $\alpha\in\NN^p_\delta$.

Up to extending $p\rest\delta$ further, we assume that $p\rest\delta$ decides the width $\width(s)$ of $s$ as well as for every $\alpha\in \NN^p_\delta$ such that $s\notin M^\delta_\alpha$, $p\rest\delta$ decides the node
\[
s_\alpha:=\text{the unique exit node from }M^{\delta}_{\alpha}\text{ below }s,
\]
as well as the widths $\width(s_\alpha)$.
Note that it is possible that $s_{\alpha}=s_{\beta}$ for $\alpha<\beta$ from $\NN^p_\delta$. The difficulty in adding each node $s_\alpha$ into $\dom(f^p_\delta)$ is that we need to find an image $t_\alpha$ for it which has a label $\lab(t)$ that separates $s_\alpha$ and $t_\alpha$ and can be added to $\EE^p_\delta$.
Enumerate the set $\NN^p_\delta$ in an increasing order as 
\[
    \alpha_0<\dots<\alpha_{n-1}.
\] 
Let $\alpha_n:=\kappa$. Note that by assumption $p\rest\delta$ has a residue to each $M^\delta_{\alpha_k}$, because $\NN^p_\delta\subseteq\EE^{p\rest\delta}_\delta$. Let $r_n:=p\rest\delta$ and by reverse recursion on $k<n$, choose $r_{k}$ such that
\[
    r_{k}\text{ is a residue of }r_{k+1}\text{ into }M^\delta_{\alpha_{k}}.
\]
We will work with this sequence of conditions $(r_k)_{k\leq n}$. By (non-reverse) recursion on $k\leq n$, we build functions $f_k$ and conditions $v_k\in\P_\delta\cap M^\delta_{\alpha_{k}}$.

\vv 

\noindent\textbf{Step $\bs{0}$:} Let $f_0:=f^p_\delta$ and $v_0:=r_0$.

\vv

\noindent\textbf{Step $\bs{k+1}$:} Assume that $f_k$ and $v_{k}\in\P_\delta\cap M^\delta_{\alpha_{k}}$ were defined and satisfy $v_{k}\leq r_{k}$. We define $f_{k+1}$ and $v_{k+1}$.
By assumption $v_{k}\in\P_\delta\cap M^\delta_{\alpha_{k}}$ and $v_{k}\leq r_{k}$. Note that also $r_k,r_{k+1},\P_\delta,M^\delta_{\alpha_k}\in M^\delta_{\alpha_{k+1}}$. Thus
\[
M^\delta_{\alpha_{k+1}}\models``r_{k}\text{ is a residue of }r_{k+1}\text{ into }M^\delta_{\alpha_k}".
\]
Hence, we may choose $v'\in\P_\delta\cap M^\delta_{\alpha_{k+1}}$ such that
\[
v'\leq v_k,r_{k+1}.
\]
If $s_{\alpha_k}\notin M^\delta_{\alpha_{k+1}}$, then let 
\[
f_{k+1}:=f_k\text{ and }v_{k+1}:=v'
\] 
and move to step $k+2$. Otherwise $s_{\alpha_k}\in M^\delta_{\alpha_{k+1}}$. We will extend $v'$ further and find $f_{k+1}\supseteq f_k$ such that
\[
\dom(f_{k+1})=\dom(f_k)\cup\{s_{\alpha_k}\}.
\]
Since $v'$ is an element in $M^\delta_{\alpha_{k+1}}$ and $\alpha_{k+1}$ is a limit point of $\EE_\delta$, there must be $\beta_k\in\EE_\delta\cap M^\delta_{\alpha_{k+1}}$ such that
\[
v'\in M^\delta_{\beta_k}.
\]
This $\beta_k$ will be the label of the image of $s_{\alpha_k}$, if $s_{\alpha_k}\in M^\delta_{\alpha_{k+1}}$.
Now $\alpha_k,\alpha_{k+1}\in\lim\EE_\delta$, $\beta_k\in\EE_\delta$ and
\[
\alpha_k<\beta_k<\alpha_{k+1}.
\]
By Lemma \ref{lem:addingmodels} applied inside the model $M^\delta_{\alpha_{i+1}}$, $v''$ defined by
\[
v''(\gamma)=\begin{cases}
    (f^{v'}_\gamma,\NN^{v'}_\gamma\cup\{\beta_k\}),\quad&\text{if }\gamma\in M^\delta_{\beta_k},\\
    (\emptyset,\emptyset),&\text{otherwise}
\end{cases}
\]
is a condition in $\P_\delta$, in particular satisfying $\beta_k\in\EE^{v''}_\delta$. Furthermore, $v''\in\P_\delta\cap M^\delta_{\alpha_{k+1}}$ and
\[
v''\leq v'\leq r_{k+1}.
\]
Consider the node $s_{\alpha_k}$. It has width $\width(s_{\alpha_k})\leq\alpha_k$ because it is an exit node from $M^\delta_{\alpha_k}$. Choose a node $t_{\alpha_k}\in T$ such that
$\height(t_{\alpha_k})=\height(s_{\alpha_k})$, $\width(t_{\alpha_k})=\alpha_k\text{ and }\lab(t)=\beta_k$.
By elementarity we may assume $t_{\alpha_k}\in M^\delta_{\alpha_{k+1}}$.
Then $t_{\alpha_k}$ is forced to be an exit node from $M^\delta_{\alpha_k}$. Furthermore, the condition $v''$ does not decide anything about $t_{\alpha_k}$, because $t_{\alpha_k}$ was chosen outside of $M^\delta_{\beta_k}$, so by extending the collapse part $f^{v''}_0$, we obtain a condition $v_{k+1}\leq v''$, still in $\P_\delta\cap M^\delta_{\alpha_{k+1}}$, that forces for every $j\leq k$, $t_{s'\wedge s}\leq_{\dot T}t_{\alpha_{j}}<_{\dot T}t_{\alpha_k}$, and for every $s'\in\dom(f^p_\delta)\cap M^\delta_{\alpha_{k+1}}$
\[
f^p_\delta(s')\wedge t_{\alpha_{k}}=t_{s'\wedge s}.
\]
Define
\[
f_{k+1}:=f_k\cup\{(s_{\alpha_k},t_{\alpha_k})\}.
\]
If $s\in M^\delta_{\alpha_{k+1}}$, then choose also a node $t_s\in (T\cap M^\delta_{\alpha_{k+1}})-\ran(f^p_\delta)$ such that $\height(t_s)=\height(s)$ and $t_{\alpha_0},\dots,t_{\alpha_{n-1}}\in\width(t_s)\times\height(t_s)$. We may assume that it is fresh in the sense that $v_{k+1}$ does not decide anything about it. Up to extending $v_{k+1}$ once more at the collapse coordinate, we may assume that it forces
\[
t_{\alpha_k}\leq_{\dot T}t_s.
\]
In this case, extend $f_{k+1}$ also by the pair $(s,t_s)$.

\noindent This ends the recursion. We have $v_n\leq p\rest\delta$ and 
\[
\hspace{1cm}v_{n}\Vdash ``f_{n}\text{ is level- and meet-preserving tree-embedding}".
\]
Now $f_n\supseteq f^p_\delta$ and $s\in\dom(f_n)$. Moreover,
\[
q:={v_n}^\smallfrown(f_n,\NN^p_\delta)
\]
is a condition that extends $p$, as wanted.
This concludes the proof of Lemma \ref{lem:dens}.
\end{proof}

\subsection{Strong Properness}\spa

\vv

We collect preliminary definitions and lemmas before presenting the proof of strong properness. As was mentioned above, the strong properness of $\P_\delta$ with respect to a model $M^\delta_\alpha$ is shown in two steps:
\begin{enumerate}
    \item if $p\in\P_\delta\cap M^\delta_\alpha$, then there is $q\leq p$ such that $\alpha\in\EE^q_\delta$,
    \item if $\alpha\in\EE^q_\delta$, then $q$ has a residue into $M^\delta_\alpha$.
\end{enumerate}
In particular, if $\alpha\in\EE^q_\delta$, then $q$ is strongly $(\P_\delta,M^\delta_\alpha)$-generic. This will imply that each poset $\P_\delta$, $\delta<\kappa^+$, is strongly proper with respect to the models
\[
(M^\delta_\alpha:\alpha\in\EE_\delta).
\]
The first step is easy. The second step is hard.

\begin{lem}\label{lem:addingmodels} Let $\delta<\kappa^+$ and $\alpha\in\EE_\delta$.
    If $p\in\P_\delta\cap M^\delta_\alpha$, then the condition $q\leq p$ defined by
    \[
    q(\gamma):=\begin{cases}
        (f^p_\gamma,\NN^p_\gamma\cup\{\alpha\}),\quad &\mbox{if }\gamma\in\delta\cap M^\delta_\alpha,\\
        (\emptyset,\emptyset) &\mbox{otherwise.}
    \end{cases}
    \] 
    is an an extension of $p$ such that $\alpha\in\EE^q_\gamma$.
\end{lem}
\begin{proof}
    It is clear that $\alpha\in\EE^q_\delta$. The fact that $q$ extends $p$ follows from the fact that the support of $p$ is contained in the model $M^\delta_\alpha$. It is straightforward to verify that $q$ is a condition: all of the clauses of the definition of the poset follow from the fact that for every $\gamma\in\delta\cap M^\delta_\alpha$, the function $f^q_\gamma$ belongs to the model $M^\gamma_\alpha$. And since each function $f^q_\gamma$ belongs to $M^\gamma_\alpha$, there are no exit nodes from $M^\gamma_\alpha$ in $\dom(f^q_\gamma)$, so the clauses in the definition of the poset $\P_\delta$ regarding exit nodes can be ignored.
\end{proof}

That was the first step. The rest of the section is devoted to the second step: building residues for conditions.

\begin{de}\label{de:trace}
    Let $\delta<\kappa^+$. The \textbf{trace} of a condition $p\in\P_\delta$ to a model $M$ is defined to be the function $[p]_M$ on $\delta$ such that
    \[
    [p]_M(\gamma)=\begin{cases}
        (f^p_\gamma\cap M,\NN^p_\gamma\cap M)\quad &\text{if }\gamma\in\delta\cap M,\\
        \emptyset &\text{if }\gamma\in\delta-M.
    \end{cases}
    \]
\end{de}
The trace of a condition is not necessarily a condition. Residues of conditions will be extensions of their traces. We write $q\leq[p]_M$ if $f^q_\gamma\supseteq f^p_\gamma\cap M\mbox{ and }\NN^q_\gamma\supseteq\NN^p_\gamma\cap M$ for all $\gamma\in\delta\cap M$, even if the trace $[p]_M$ was not a condition.

\begin{notation}
    For a finite set $E\subseteq\kappa^+\times\kappa$ of pairs of ordinals we write
    \begin{align*}
        E_\gamma:=\{\beta:(\gamma,\beta)\in E\}
    \end{align*}
    and when the set $E$ is clear from context, we also denote by $\bs{(\gamma,\beta^+)}$ the pair where
    \begin{align*}
        \beta^+:=\begin{cases}
            \text{the least element of }E_\gamma\text{ strictly above }\beta,\quad & \text{if }\beta<\max(E_\gamma),\\
            \kappa, &\text{if }\beta=\max(E_\gamma),
            \end{cases}
    \end{align*}
    for any $\beta<\kappa$.
    Furthermore, we denote
    \[
    E\rest\delta:=E\cap(\delta\times\kappa).
    \]
    
\end{notation}

Recall from item (\ref{prop:1tree:free}) of Proposition \ref{prop:1tree} that a condition $f\in\Coll(\omega,<\kappa)$ is \textbf{nice with respect to a node $\bs{t}$} if for any $\alpha$ with $t\notin V_\alpha$, any $w\in \Coll(\omega,<\kappa)\cap V_\alpha$ and any $\bar{t}\in\width(t)\times\height(t)$, if $w\leq p\cap V_\alpha$ and $w\Vdash ``\pi^p(t)<_{\dot T}\bar{t}"$, then there is $q\leq p$ such that $q\cap v_\alpha=p\cap V_\alpha$, $q$ is minimal in the sense that $\pi^q(t')=\pi^p(t')$ for any $t'\neq t$, and that any common extension of $w$ and $q$ forces $\bar{t}<_{\dot T}t$. Recall also that the set of conditions that are nice with respect to $t$ for every $t\in A$ for any $A\subseteq T$, are dense in $\Coll(\omega,<\kappa)$.

\begin{de}
    Let $\delta<\kappa^+$, $p\in\P_\delta$ and $\alpha\in\EE^p_\delta$. A \textbf{residue system for $\bs{p}$ into $\bs{M^\delta_\alpha}$} is a tuple
    \[
    \vec{r}_E=(r_{(\gamma,\beta)}:(\gamma,\beta)\in E)
    \]
    indexed by a finite set $E\subseteq\kappa^+\times\kappa$ such that if we denote 
    \[
        r_{(\gamma,\kappa)}:=p\rest\gamma,
    \] 
    then the following are satisfied:
    \begin{enumerate}
        \item $r_{(\gamma,\beta)}\in\P_\gamma\cap M^\gamma_\beta$ for every $(\gamma,\beta)\in E$.
        \item\label{item:ressys:3} $E$ is the least set such that $E_\delta=\{\alpha\}$ and $E$ is \textbf{closed} in the following sense: for every $(\gamma,\beta)\in E$ and $\xi\in \sp(r_{(\gamma,\beta^+)})\cap M^\gamma_\beta$ with $\xi\geq 1$, 
        \[
        \EE^{r_{(\gamma,\beta^+)}}_\xi-\beta\subseteq E_\xi.
        \]
        Furthermore, if $\beta^+\neq\kappa$, then $\beta^+\in E_\xi$.
        \item\label{item:ressys:1} 
        Let $(\gamma,\beta),(\gamma',\beta')\in E$. 
        If $\beta<\beta'$, or if $\beta=\beta'$ and $\gamma>\gamma'$, then
        \[
        r_{(\gamma,\beta)}\rest\min\{\gamma,\gamma'\}\leq[r_{(\gamma',\beta')}\rest\min\{\gamma,\gamma'\}]_{M^\gamma_\beta}
        \]
        In particular, $\bigcup_{x\in E}f^{r_x}_0$ is a condition in $\Coll(\omega,<\kappa)$.
        \item\label{item:ressys:5} For every $(\gamma,\beta)\in E$ and every non-negative  
        $\xi\in \gamma\cap M^\gamma_\beta$
        and for every pair $(s,t)\in f^{r_{(\gamma,\beta^+)}}_{\xi}$ of exit nodes from $M^\xi_\beta$, if $\rho$ is the maximal ordinal in $E_\xi$ such that $t\notin M^\xi_\rho$, then $r_{(\xi,\rho)}$ decides the ``implicit preimage" for the projection of the image of $s$, i.e. the node
        \[
        \bar{s}:=\text{the predecessor of }s\text{ at the height of }\pi^E(t),
        \]
        where $\pi^E(t):=\pi^{f}(t)$ for $f=f^p_0\cup\bigcup_{x\in E}f^{r_x}_0$. 
        \item\label{item:ressys:6} For every $(\gamma,\beta)\in E$ and $\xi\in\gamma\cap M^\gamma_\beta$, the collapse part $f^{r_{(\gamma,\beta^+)}}_0$ is nice with respect to $t$ for every $t\in\ran(f^{r_{(\gamma,\beta^+)}}_\xi)$ that is an exit node from $M^\gamma_\beta$.
    \end{enumerate}
    We call the condition $r_{(\delta,\alpha)}$ the \textbf{root condition} of the system $\vec{r}_E$.
\end{de}

\begin{notation}
    If $\vec{r}_E$ is a residue system for $p\in\P_\delta$ into $M^\delta_\alpha$ and there is no danger of confusion of $p$, we denote
    \[
    r_{(\gamma,\kappa)}:=p\rest\gamma,
    \]
    for every $\gamma<\delta$.
\end{notation}

\begin{lem}
    If $\vec{r}_E$ is a residue system for $p$ into $M^\delta_\alpha$, then $f^p_0\cup \bigcup_{x\in E}f^{r_x}_0$ is a condition in $\Coll(\omega,<\kappa)$.
\end{lem}
\begin{proof}
    Follows from item (\ref{item:ressys:1}) of definition of residue system.
\end{proof}

\begin{lem}
    If $\vec{r}_E$ is a residue system for $p$ into $M^\delta_\alpha$ and $r\in\P_\delta\cap M^\delta_\alpha$ extends $r_{(\delta,\alpha)}$, then the system obtained from $\vec{r}_E$ by replacing $r_{(\delta,\alpha)}$ by $r$ is a residue system for $p$ into $M^\delta_\alpha$.
\end{lem}

\begin{de}
    A \textbf{path} is a finite sequence $(\gamma_0,\beta_0),\dots,(\gamma_n,\beta_n)$ with
\begin{enumerate}
    \item $\gamma_{k+1}\in\gamma_k\cap M^{\gamma_k}_{\beta_k}$,
    \item $\beta_{k+1}\in\EE_{\gamma_{k+1}}-\beta_k$.
\end{enumerate}
A \textbf{path from $\bs{M^\delta_\alpha}$ to $\bs{M^\gamma_\beta}$} is a path with $(\gamma_0,\beta_0)=(\delta,\alpha)$ and $(\gamma_n,\beta_n)=(\gamma,\beta)$. If $E$ is a set of pairs, then a \textbf{path in $\bs{E}$} is a path whose each member belogns to $E$.
\end{de}

\begin{lem}
    If $\vec{r}_E$ is a residue system for a condition $p\in\P_{\delta}$ into $M^{\delta}_{\alpha}$ and $(\gamma,\beta)\in E$, then there is a path from $M^{\delta}_{\alpha}$ to $M^\gamma_\beta$.
\end{lem}

\begin{lem}\label{lem:claim52}\spa 

\begin{enumerate}
    \item\label{item1:lem:paths} If there is a path from $M^\delta_\alpha$ to $M^\gamma_\beta$, then $\gamma\cap M^{\delta}_{{\alpha}}\subseteq M^\gamma_\beta$. 
    \item\label{item2:lem:paths} If there is a path from $M^\delta_\alpha$ to $M^\gamma_\beta$ and and $\beta\leq \alpha^+$, then $\gamma\in M^\delta_{\alpha^+}$, where $\alpha^+$ is the next element of $\EE_\delta$ strictly above $\alpha$.
\end{enumerate}
    
\end{lem}
\begin{proof}
    We prove item (\ref{item1:lem:paths}). The proof is by induction on the length of the path. To this end, fix a path
    \[
    (\gamma_0,\beta_0),\dots,(\gamma_{n+1},\beta_{n+1}).
    \]
    with $(\gamma_0,\beta_0)=(\delta,\alpha)$.
    We claim that 
    \[
    \gamma_{n+1}\cap M^{\delta}_{\alpha}\subseteq M^{\gamma_{n+1}}_{\beta_{n+1}}.
    \]
    Let $\xi\in\gamma_{n+1}\cap M^{\delta}_{\alpha}$. Let $\psi$ be the $<_\theta$-least bijection from $\kappa$ to $\gamma_{n+1}$. Then $\psi\in M^{\gamma_{n+1}}_{\beta_{n+1}}$ but also $\psi\in M^{\gamma_{n}}_{\beta_{n}}$ since $\gamma_{n+1}\in M^{\gamma_n}_{\beta_n}$. By the induction hypothesis we have $\gamma_n\cap M^{{\delta}}_{\alpha}\subseteq M^{\gamma_n}_{\beta_n}$, so $\xi\in M^{\gamma_n}_{\beta_n}$ too. And thus $\bar{\xi}:=\psi^{-1}(\xi)\in\beta_n\subseteq\beta_{n+1}$. And then $\xi=\psi(\bar{\xi})\in M^{\gamma_{n+1}}_{\beta_{n+1}}$. The fact that $\P_\gamma\cap M^{\delta}_{\alpha}\subseteq \P_\gamma\cap M^\gamma_\beta$ follows immediately.

    The proof of item (\ref{item2:lem:paths}) is proved the same way by induction on the length of the path, using appropriate injection $\psi$.
\end{proof}

\begin{lem}\label{lem:E''}
    If $\vec{r}_E$ is a residue system for $p$ into $M^\delta_\alpha$ and $\gamma=\max(\sp(p)\cap M^\delta_\alpha)$, then $E-\{(\delta,\alpha)\}$ can be written as the disjoint union of the sets
    \[
    E^\beta:=E\cap(M^\gamma_{\beta^+}-M^\gamma_\beta)
    \]
    where $\beta\in E_\gamma$, and furthermore, $\vec{r}_{E^\beta}$ is a residue system in $M^\gamma_{\beta^+}$ for $r_{(\gamma,\beta^+)}$ into $M^\gamma_\beta$.
\end{lem}
\begin{proof}
    Enumerate the set $E_\gamma$ as $\alpha=\beta_0<\dots<\beta_{n-1}$ and let $\beta_n:=\kappa$ and denote $E^k:=E\cap(M^\gamma_{\beta_{k+1}}-M^\gamma_{\beta_k})$. Fix some $k\leq n$. Then $\vec{r}_{E^k}\in M^\gamma_{\beta_{k+1}}$, since if $(\xi,\rho)\in E^k$, then $r_{(\xi,\rho)}\in M^\xi_\rho\subseteq M^\gamma_{\beta_{k+1}}$. We claim that $\vec{r}_{E^k}$ is a residue system for $r_{k+1}:=r_{(\gamma,\beta_{k+1})}$ into $M^\gamma_{\beta_k}$.

    We say that a path $(\gamma_k,\beta_k)_{k\leq m}$ is a \textbf{path in $\bs{\vec{r}_E}$} if it is a path in $E$ that satisfies additionally: $\gamma_{k+1}\in\sp(r_{(\gamma_k,\beta_k^+)})\cap M^{\gamma_k}_{\beta_k}$ and $\beta_{k+1}\in \EE^{r_{(\gamma_k,\beta^+_k)}}_{\gamma_{k+1}}-\beta_k$ for every $k\leq m$.

    \begin{claim}\label{claim:E^k}
        If $(\xi,\rho)\in E^k$, then $r_{(\xi,\rho^+)}\in M^\gamma_{\beta_{k+1}}$.
    \end{claim}
    \begin{proof}[Proof of Claim \ref{claim:E^k}]
        The claim is vacuously true if $\xi=\gamma$, since then $\rho=\beta_k$ and $r_{(\gamma,\beta_k)}\in M^\gamma_{\beta_k}\subseteq M^\gamma_{\beta_k}$. 
        We show by induction on the length of the path that if $(\xi_i,\rho_i)_{i\leq m}$ is a path in $\vec{r}_E$ with $\rho_m<\beta_{k+1}$, then $\rho_m^+\leq\beta_{k+1}$. Fix a path $(\xi_i,\rho_i)_{i\leq m+1}$ in $\vec{r}_E$ such that $\rho_{m+1}<\beta_{k+1}$. We claim that $\rho_{m+1}^+\leq\beta_{k+1}$. Now $\rho_m\leq \rho_{m+1}$, so by the induction hypothesis $\rho_m^+\leq\beta_{k+1}$. In paticular $\rho^+_m<\kappa$, so
        \[
        \rho^+_m\in E_{\xi_{m+1}}.
        \]
        But note that $\rho_{m+1}\in(\EE^{r_{(\xi_m,\rho^+_m)}}_{\xi_{m+1}}-\rho_m)\subseteq M^{\xi_m}_{\rho^+_m}$, so $\rho_{m+1}<\rho^+_m$. This implies that
        \[
        \rho^+_{m+1}\leq\rho^+_m\leq\beta_{k+1}.
        \]
        Thus, if $(\xi,\rho)\in E^k$, then there is path from $M^\delta_\alpha$ to $M^\xi_\rho$ in $\vec{r}_E$, and so $\rho<\beta_{k+1}$ implies $\rho^+\leq\beta_{k+1}$, which in turn implies $r_{(\xi,\rho^+)}\in M^\gamma_{\beta_{k+1}}$.
    \end{proof}

    Now, the above claim implies that $\vec{r}_{E^k}$ satisfies item (\ref{item:ressys:3}) of the definition of residue system, namely, the following: $E^k$ is the least set such that $E^k_\gamma=\{\beta_k\}$ and it is closed for the following: if $(\xi,\rho)\in E^k$ and $\xi'\in\sp(r_{(\xi,\rho^+)})\cap M^\xi_\rho$, then
    \[
    \EE^{r_{(\xi,\rho^+)}}_\xi-\rho\subseteq E^k_\xi,
    \]
    and if $\rho^+<\kappa$, then also $\rho^+\in E^k_\xi$. It also also implies that if $(\xi,\rho)\in E^k$, then there is a path in $\vec{r}_{E^k}$ from $M^\gamma_{\beta_k}$ to $M^\xi_\rho$, and furthermore, that if $(\xi,\rho)\in E^{k+1}-E^k$, then $\rho\geq\beta_{k+1}$.

    It remains to be shown that items (\ref{item:ressys:5}) and (\ref{item:ressys:6}) are satisfied. 

    \begin{claim}\label{claim:proj}
        If $t\in T\cap(M^\gamma_{\beta_{k+1}}-M^\gamma_{\beta_k})$, then $\pi^E(t)=\pi^{E^k}(t)$.
    \end{claim} 
    \begin{proof}[Proof of Claim \ref{claim:proj}]
        Follows from the fact that the branch below $t$ is decided by the poset $\Coll(\omega,<\kappa)\cap(M^\gamma_{\beta_{k+1}}-M^\gamma_{\beta_k})$, and from the fact that if $(\xi,\rho)\in E^{k+1}-E^k$, then $\rho>\beta_{k+1}$, and thus for any $(\xi',\rho')\in E^k$, the collapse conditions extend: $f^{r_{(\xi',\rho')}}_0\leq [f^{r_{(\xi,\rho)}}_0]_{V_{\rho'}}$.
    \end{proof}
    
    Claim \ref{claim:proj} implies that items (\ref{item:ressys:5}) and (\ref{item:ressys:6}) are satisfied. This concludes the proof of the lemma.

\end{proof}

In general, if $\vec{r}_E$ is a residue system for $p$ into $M^\delta_\alpha$ and $t\in T$ is a node, we denote by $\pi^E(t)$ the node projection $\pi^f(t)$ where
\[
f=f^p_0\cup\bigcup_{x\in E}f^{r_x}_0.
\]

\begin{lem}[Flexibility Lemma]\label{lem:freemodifs} Let $\delta<\kappa$, $p\in\P_\delta$ and $\alpha\in\EE^p_\delta$. Assume that $\vec{r}_E$ is a residue system for $p$ into $M^\delta_\alpha$. Suppose that $t\in T$ is a node with $\lab(t)\geq\alpha$. If $w\in\P_\delta\cap M^\delta_\alpha$ extends $r_{(\delta,\alpha)}$ and $\bar{t}\in\width(t)\times\height(t)$ is a node and \[
w\Vdash ``\pi^E(t)<_{\dot T}\bar{t}",
\]
then there is a condition $p'\leq p$ and residue system $\vec{r}'_E$ for $p'$ into $M^\delta_\alpha$ with
\begin{enumerate}
    \item the root conditions are equal, $r'_{(\delta,\alpha)}=r_{(\delta,\alpha)}$,
    \item any common extension of $w$ and $p'$ forces $``\bar{t}<_{\dot T}t"$.
\end{enumerate}
    
\end{lem}
\begin{proof}
    Let $g\leq f^p_0$ be as in item (\ref{prop:1tree:free}) from Proposition \ref{prop:1tree}. Obtain $p'$ and $\vec{r}'_{E}$ from $p$ and $\vec{r}_E$, respectively, by extending each condition at their collapse coordinate by letting
    \[
    f^{r'_{(\gamma,\beta)}}_0:=f^{r_{(\gamma,\beta)}}_0\cup (g\cap M^\gamma_\beta),
    \]
    for every $(\gamma,\beta)\in E$. 
    We claim that $\vec{r'}_E$ is a residue system for $p'$ into $M^\delta_\alpha$. We verify items (\ref{item:ressys:5}) and (\ref{item:ressys:6}). We show first that $t$ cannot be in the range of any $f^{r_{(\gamma,\beta^+)}}_\xi$ and exit node from $M^\gamma_\beta$, for $(\gamma,\beta)\in E$. For if it was, the preimage of $t$ would also be an exit node from $M^\gamma_\beta$. Since $\beta\in\NN^{r_{(\gamma,\beta^+)}}_\xi$, this would mean that the preimage of $t$ should also be an element in $V_{\lab{t}}$, by definition of the poset, which is impossible since $\lab(t)\leq\alpha\leq\beta$. Thus, items (\ref{item:ressys:5}) and (\ref{item:ressys:6}) cannot fail by minimality of the extension $g\leq f^p_0$.
\end{proof}

\begin{lem}\label{lem:spsysisres}
    Let $\delta<\kappa^+$. If $\vec{r}_E$ is a residue system for $p$ into $M^\delta_\alpha$, then for every $(\gamma,\beta)\in E$, the condition $r_{(\gamma,\beta)}$ is a residue for $r_{(\gamma,\beta^+)}$ into $M^\gamma_\beta$.
\end{lem}
\begin{proof}
The proof is by induction on $\delta$. By the induction hypothesis and Lemma \ref{lem:E''}, it is enough to show that the root condition $r_{(\delta,\alpha)}$ is a residue for $p$, since $r_{(\delta,\alpha^+)}=p$. Let $w\in\P_\delta\cap M^\delta_\alpha$ extend $r_{(\delta,\alpha)}$. We find a common extension of $w$ and $p$.

\vv

\noindent \textbf{Base case and limit case:}

\vv

The case $\delta=1$ follows from the fact that the working parts of conditions in $\P_1$ are in $\Coll(\omega,<\kappa)$, so the pointwise union $w\cup p$ is a common extension of $w$ and $p$. See \ref{ex:coll}. The case when $\delta$ is a limit ordinal follows from the induction hypothesis: Let $\gamma<\delta$ be such that $\sp(w)\subseteq\gamma$. We can choose such $\gamma$ in $M^\delta_\alpha$. By Lemma \ref{lem:structureofmodels}, $\gamma\cap M^\delta_\alpha\subseteq M^\gamma_\alpha$. Thus $w\rest\gamma\in M^\gamma_\alpha$. Define the set $E'$ from $E$ by changing the pair $(\delta,\alpha)$ to $(\gamma,\alpha)$, and define the system $\vec{r}_{E'}$ from $\vec{r}_E$ by replacing the root condition $r_{(\delta,\alpha)}$ by the condition $r_{(\gamma,\alpha)}:=r_{(\delta,\alpha)}\rest\gamma$. The resulting system is a residue system for $p\rest\gamma$ into $M^\gamma_\alpha$. Now $w\rest\gamma$ extends the root condition $r_{(\gamma,\alpha)}$ so  by the induction hypothesis, $w\rest\gamma$ is compatible with $p\rest\gamma$. Let $q\leq w\rest\gamma,p\rest\gamma$. Then $q^\smallfrown p\rest[\gamma,\delta)$ is a common extension of $w$ and $p$.

\vv 

\noindent\textbf{Successor case $\boldsymbol{\delta+1}$:} 

\vv

Now $w,p\in\P_{\delta+1}$. We will find a condition $v$ in $\P_\delta$ that is a common extension of $ w\rest\delta$ and $p\rest\delta$, and a function $f$ extending the functions $f^w_\delta$ and $f^p_\delta$, in such a way that $q:=v^\smallfrown(f,\NN^w_\delta\cup\NN^p_\delta)$ will be the desired common extension. For instance, this $v$ must decide the meets in the set $\dom(f^w_\delta)\cup\dom(f^p_\delta)$, and the meets must be contained in the domain of $f$. The condition $v$ and the function $f$ will be defined recursively, by climbing up the models $M^\delta_{\beta_k}$. It is important to do this gradually, recursively, in order to avoid problems discussed in the introduction of this section. The main idea is the following: if we decide meets of nodes in $\dom(f^w_\delta)$ and a node $s\in\dom(f^p_\delta)-M^\delta_\alpha$ in a model $M^\delta_{\beta}$ that contains $s$ but not its image $f^p_\delta(s)$, we do not touch the node projection $\pi(f^p_\delta(s))$, and are later free to modify the branch below $f^p_\delta(s)$ to correspond to the picture of the branch below $s$ in the domain-side.

\vv

By minimality of $E$, we have $E_\delta=\EE^p_\delta-\alpha$. Enumerate $E_\delta$ in a strictly increasing order as
\[
\alpha=\beta_0<\dots<\beta_{n-1}.
\]
Denote $\beta_{n}:=\kappa$, $r_{n}:=p\rest\delta$ and for each $k< n$, abbreviate
\begin{align*}
    & r_k:=r_{(\delta,\beta_k)}\\
    & E_k:=E\cap(M^\delta_{\beta_{k+1}}-M^\delta_{\beta_k})\cap(\delta+1\times\kappa).
\end{align*}
We may assume that we added the pair $(\delta,\alpha)$ to $E_0$ and denote $r_0:=r_{(\delta+1,\alpha)}\rest\delta$. By Lemma \ref{lem:E''}, the system $\vec{r}_{E_k}$ is a residue system for $r_{k+1}$ into $M^\delta_{\beta_k}$, and thus by the induction hypothesis, $r_k$ is a residue of $r_{k+1}$ into $M^\delta_{\beta_k}$.
Also, for $k\leq n$, let
\[
    X_{k}:=\{s\in\dom(f^p_\delta):s\text{ is exit node from }M^\delta_\alpha\text{ and }f^p_\delta(s)\in M^\delta_{\beta_{k+1}}-M^\delta_{\beta_{k}}\}.
\]
Then every node in $\dom(f^p_\delta)$ that is an exit node from $M^\delta_\alpha$ belongs to some $X_k$, $k< n$, and $X_k\subseteq M^\delta_{\beta_k}$. The fact that $X_k\subseteq M^\delta_{\beta_k}$ follows from the fact that if $s\in X_k$, then $\lab(f^p_\delta(s))\leq\beta_k$ and thus $s\in V_{\lab(t)}\subseteq M^\delta_{\beta_k}$.

We will define by recursion on $k\leq n$ conditions $v_k$ and functions $f_k$.    
\vv
        
\noindent\textbf{Step $\bs{0}$:} Let $v_0:=w\rest\delta$ and $f_0:=f^w_\delta$.

\vv

\noindent\textbf{Step $\bs{k+1}$:} Assume that $v_k$ and $f_k$ are defined. Suppose that they satisfy:
\begin{enumerate}
    \item $v_k\in\P_\delta\cap M^\delta_{\beta_k}$ and $v_k\leq r_k$.
    \item $f_k$ is a function in $M^\delta_{\beta_k}$ such that
    \[
    {v_k}^\smallfrown(f_k,\NN^w_\delta)\in\P_{\delta+1}
    \]
    and the domain of $f_k$ is the closure of the set $\dom(f^w_\delta)\cup\bigcup_{j<k}X_j$ under meets and under taking exit nodes from models $M^\delta_{\alpha'}$ where $\alpha'\in\NN^w_\delta$, as decided by $v_k$.
\end{enumerate}
Up to extending $v_k$ inside the model $M^\delta_{\beta_k}$, we may assume that it decides meets $s\wedge s'$ where $s\in X_k$ and $s'\in\dom(f^w_\delta)$, as well as their \textbf{implicit images}:
\[
\bs{t_{s\wedge s'}}:=\text{the node below }f^w_\delta(s')\text{ at the height of }s\wedge s',
\]
as well as for every $\alpha'\in\NN^w_\delta$ and $s\in X_k$, the node
\[
\bs{\bar{s}_{\alpha'}}:=\text{the node }s'\leq s\text{ which is an exit node from }M^\delta_{\alpha'}.
\]
Furthermore, up to extending $r_k$ even further, by Lemma \ref{lem:dens}, we may assume that for each $\bar{s}_{\alpha'}$ such that $\bar{s}_{\alpha'}\neq s$, there is a node $t_{\bar{s}_{\alpha'}}$ such that if we let $f'_k$ to be the function $f_k\cup\{(\bar{s}_{\alpha'},\bar{t}_{\alpha'}):\alpha'\in\NN^w_\delta,\bar{s}_{\alpha'}\neq s, s\in X_k\}$, then
\[
{v_k}^\smallfrown(f'_k,\NN^w_\delta)\in\P_{\delta+1}.
\]
Since $\width(f^p_\delta(s))\in\NN^p_\delta$, we automatically have $t_{\bar{s}_{\alpha'}}\in\width(f^p_\delta(s))\times\width(f^p_\delta(s))$ for every $s\in X_k$. Note also that $r_k\Vdash ``\pi^{E_k}(f^p_\delta(s))\leq f_k(\tilde s)"$ if $\tilde s\notin\dom(f^{r_{(\delta+1,\alpha)}}_\delta)$, for any $\tilde s\in\{s\wedge s':s'\in\dom(f^w_\delta)\}\cup(\{\bar{s}_{\alpha'}:\alpha'\in\NN^w_\delta\}\cap V_\alpha)$.

The goal is to define $v_{k+1}$ and $f_{k+1}$. We will work in the model $M^\delta_{\beta_{k+1}}$. We will first look at nodes in the set $f^p_\delta[X_k]$ and modify branches below them, as follows.

Recall that $\vec{r}_{E_k}$ is a residue system for $r_{k+1}$ into $M^\delta_{\beta_k}$ with root condition $r_k$. By the Flexibility Lemma \ref{lem:freemodifs}, since $\lab(t)\leq\beta_k$ for every $t\in f^p_\delta[X_k]$, we find a condition $r'_{k+1}\leq r_{k+1}$ and a residue system $\vec{r}'_{E_k}\leq \vec{r}_{E_k}$ such that:
\begin{enumerate}
    \item the root conditions are equal $r'_{\delta,\beta_k}=r_k$,
    \item for any $s\in X_k$, $s'\in\dom(f^w_\delta)$ and $\alpha'\in\NN^w_\delta$, any common extension of $v_k$ and $r'_{k+1}$ forces the following:
    \begin{enumerate}
        \item $f^p_\delta(s)\wedge f^w_\delta(s')=t_{s\wedge s'}$,
        \item $t_{\bar{s}_{\alpha'}}\leq _{\dot T}f^p_\delta(s)$.
    \end{enumerate} 
\end{enumerate}
By the induction hypothesis $r_k$ is a residue for $r'_{k+1}$ and since $v_k\leq r_k$, it must be compatible with $r'_{k+1}$. Let $v_{k+1}\leq v_k,\tilde r$. Applying the induction hypothesis inside the model $M^\delta_{\beta_{k+1}}$, we may assume that $v_{k+1}\in\P_\delta\cap M^\delta_{\beta_{k+1}}$. 
Define then $f_{k+1}$ to be the function consisting of the following sets:
\begin{itemize}
    \item $f_k$,
    \item $\{(s\wedge s',t_{s\wedge s'}):s\in X_k\text{ and }s'\in\dom(f^w_\delta)\}$,
    \item $\{(\bar{s}_{\alpha'},t_{\bar{s}_{\alpha'}}):s\in X_k, \alpha'\in\NN^w_\delta\}$.
\end{itemize}
Then $f_{k+1}$ is an injective function and by construction $v_{k+1}$ forces that it is level- and meet-preserving tree-embedding, and moreover,
\[
{v_{k+1}}^\smallfrown(f_{k+1},\NN^w_\delta)\in\P_{\delta+1}.
\]
This ends the recursion. Finally, look at the condition $v_{n}$ and the function $f_{n}$. Define
\[
q:={v_{n}}^\smallfrown(f_{n}\cup f^p_\delta,\NN^w_\delta\cup\NN^p_\delta).
\]
By construction, $q$ is a condition in $\P_{\delta+1}$ that extends both $w$ and $p$. Hence $r_{(\delta+1,\alpha)}$ is a residue of $p$. This ends the proof of Lemma \ref{lem:spsysisres}.
        
\end{proof}

\begin{lem}\label{lem:spconstres}
    If $\alpha\in\EE^p_\delta$, then $p$ has a residue system into $M^\delta_\alpha$.
\end{lem}
\begin{proof}
The proof is by induction on $\delta<\kappa^+$. First of all, up to extending $p(0)$, we may assume that $f^p_0$ is nice with respect to every $t\in\bigcup_{\gamma<\delta}\ran(f^p_\gamma)$, in the sense of item (\ref{prop:1tree:free}) of Proposition \ref{prop:1tree}.

\vv

\noindent\textbf{Base case $\bs{\delta=1}$:} If $p\in\P_1$ and $\alpha\in\EE^p_\delta$, we may let the index set to be $E:=\{(1,\alpha)\}$ and the condition $r_{(1,\alpha)}$ to be the pairwise intersection $(f^p_0\cap V_\alpha,\NN^p_0\cap\alpha)$. This system is as wanted.

\vv

\noindent\textbf{Limit case $\bs{\delta}$:} Suppose that $\delta$ is a limit ordinal. Let $p\in\P_\delta$ and $\alpha\in\EE^p_\delta$. Let $\gamma\in\delta\cap M^\delta_\alpha$ be such that $\sp(p)\cap M^\delta_\alpha\subseteq\gamma$. Such $\gamma$ exists by Lemma \ref{lem:structureofmodels}. By the induction hypothesis $p\rest\gamma$ has a residue system $\vec{r}_D$ into $M^\gamma_\alpha$. We define $E$ from $D$ by replacing $(\gamma,\alpha)$ by $(\delta,\alpha)$ and letting $r_{(\delta,\alpha)}:={r_{(\gamma,\alpha)}}^\smallfrown(\emptyset,\emptyset),\dots$. Then $r_{(\delta,\alpha)}\in\P_\delta\cap M^\delta_\alpha$ and the system $\vec{r}_E$ is as wanted.

\vv 

\noindent\textbf{Successor case $\boldsymbol{\delta+1}$:} 

\vv

Let $p\in\P_{\delta+1}$ and $\alpha\in\EE^p_{\delta+1}$. We build a residue system for $p$. First, enumerate the finite set $\EE^p_\delta-\alpha$ in a strictly increasing order as
\[
\alpha=\beta_0<\beta_1<\dots<\beta_{n-1}.
\]
Abbreviate $\beta_{n}:=\kappa$. For every $k<n$, denote
\[
X_k:=\{s\in\dom(f^p_\delta):s\text{ is an exit node from }V_\alpha\text{ and }f^p_\delta(s)\in M^\delta_{\beta_{k+1}}- M^\delta_{\beta_k}\}.
\]
Then $\lab(f^p_\delta(s))\leq\beta_k$ for every $s\in X_k$ and thus $X_k\subseteq M^\delta_{\beta_k}$.
We proceed by reverse recursion on $k\leq n$.

\vv 

Let $r_{n}:=p\rest\delta$. At step $k<n$, suppose that $r_{k+1}$ has been defined and satisfies $r_{k+1}\in\P_\delta\cap M^\delta_{\beta_{k+1}}$.
By the induction hypothesis $r_{k+1}$ has a residue system $\vec{r}_{E_{k}}$ into $M^\delta_{\beta_{k}}$. By working inside the model $M^\delta_{\beta_{k+1}}$, by elementarity, we may assume that the residue system $\vec{r}_{E_k}$ is an element of $M^\delta_{\beta_{k+1}}$. Denote its root condition $r_{(\delta,\beta_k)}$ by $r_k$.
Then $r_k\in \P_\delta\cap M^\delta_{\beta_k}$. Up to extending $r_k$ inside the model $M^\delta_{\beta_k}$, we may assume that it decides, for every $s\in X_k$, the ``implicit preimage" of the node $s$, i.e. the predecessor $\bar{s}$ of $s$ at the height of \[t_{\bar{s}}:=\pi^g(f^p_\delta(s)),\] where
$g=f^{r_{k+1}}_0\cup\bigcup_{x\in E_k}f^{r_x}_0$.
Furthermore, up to extending the corodinate zero $f^{r_k}_0$, still inside the model $M^\delta_{\beta_k}$, we may assume that it is nice with respect to every $t\in X_k$ in the sense of item (\ref{prop:1tree:free}) of Proposition \ref{prop:1tree}.

This ends the recursion. We have now chosen the conditions $r_k$, $k\leq n$, and residue system $\vec{r}_{E_k}$ for $r_{k+1}$ into $M^\delta_{\beta_k}$, as well as nodes $(\bar{s},\bar{t}_s)$ for each $s\in X_k$. We look at the final condition $r_0\in\P_\delta\cap M^\delta_\alpha$ and define 
\[
f:=(f^p_\delta\cap M^\delta_\alpha)\cup\bigcup_{k\leq n}\{(\bar{s},\bar{t}_s):s\in X_k\}.
\]
Let
\[
r:={r_0}^\smallfrown(f,\NN^p_\delta\cap\alpha).
\] 
We first verify that $r$ is a condition in $\P_{\delta+1}\cap M^{\delta+1}_\alpha$.

\begin{claim}\label{claim:stp:cond}
    $r\in\P_{\delta+1}\cap M^{\delta+1}_\alpha$.
\end{claim}
\begin{proof}[Proof of Claim \ref{claim:stp:cond}]
To see $r\in M^{\delta+1}_\alpha$, it suffices to see that $r_0,f$ and $\NN^p_\delta\cap\alpha$ are in $M^{\delta+1}_\alpha$. This is clear for $r_0$ and $\NN^p_\delta\cap\alpha$. We argue that the function $f$ is in $M^{\delta+1}_\alpha$. Since $\alpha\in\NN^p_\delta$, the model $M^{\delta}_\alpha$ and thus also $M^{\delta+1}_\alpha$ must be closed under the function $f^p_\delta$, whence $f^p_\delta\cap M^{\delta}_\alpha\in M^{\delta+1}_\alpha$. And since each $s\in X_k$, for $k\leq n$, is an exit node from $M^\delta_\alpha$, so are the images $f^p_\delta(s)$, and thus the pairs $(\bar{s},\bar{t}_s)$ must belong to $M^\delta_\alpha$. Thus $f\in M^{\delta+1}_\alpha$. Hence $r\in M^{\delta+1}_\alpha$.

We then argue that $r\in\P_{\delta+1}$. Since $r\rest\delta=r_0\in\P_\delta$, it suffices to verify that it forces all the relevant items about $f$ and $\NN^p_\delta\cap\alpha$ from the definition of the poset. This follows straightforwardly using the fact that $r_0$ is a residue for $p\rest\delta$. For instance, to show that $r_0$ decides meets in the set $\dom(f)$, let $s,s'\in \dom(f)$. Note first that if $s,s'\in\dom(f^p_\delta)$, then there is automatically some $\bar{s}\in\dom(f)$ such that $p\rest\delta\Vdash``s\wedge s'=\bar{s}"$ and $r_0$ must force the same thing, for otherwise it cannot be a residue of $p\rest\delta$. On the other hand, if $s\notin\dom(f^p_\delta)$, then $s=\bar{s}_0$ for some $s_0\in X_k$, and then 
\[
r_k\Vdash`` \bar{s}_0\wedge s'=s_0\wedge s'".
\]
This implies that there is $\bar{s}\in\dom(f^p_\delta)\cap M^\delta_\alpha$ such that 
\[
p\rest\delta\Vdash``\bar{s}=s_0\wedge s'".
\]
And then $r_0$ must also force $``\bar{s}=\bar{s}_0\wedge s'"$, being a residue for both $r_k$ and $p\rest\delta$. Hence $r_0$ must decide meets in $\dom(f)$. The remaining properties are verified similarly. 

Note that the last clauses in the definition of the poset regarding models $M^\delta_\beta$, $\beta\in\NN^p_\delta\cap\alpha$ and exit nodes from them in $\dom(f)$ are vacuously satisfied due to the fact that we only added nodes below nodes in $\dom(f^p_\delta)$: for instance, if $s\in X_k$ and $\bar{s}\notin M^\delta_\beta$, for some $\beta\in\NN^p_\delta\cap\alpha$, then the unique exit node from $M^\delta_\beta$ below $\bar{s}$ is the same node as the unique exit node below $s$, which must be decided by $p\rest\delta$ and thus by $r_0$ and must be in $\dom(f^p_\delta)\cap M^\delta_\alpha$, as $p$ is a condition. This ends the proof of Claim \ref{claim:stp:cond}.

\end{proof}

We are ready to define the residue system for $p$ into $M^\delta_{\alpha+1}$. The index set $E$ will be the set $\bigcup_{k\leq n}E_k$ where the pair $(\delta,\alpha)$ is replaced by the pair $(\delta+1,\alpha)$, and for each $(\gamma,\beta)\in E_k$ such that $\beta^+<\kappa$, the pair $(\gamma,\beta_{k+1})$ is added to $E$. For such pairs, we let $r_{(\gamma,\beta_{k+1})}:=r_{k+1}\rest\gamma$. Let $r_{(\delta+1,\alpha)}:=r$. For $(\gamma,\beta)\in E_k$, we have already defined $r_{(\gamma,\beta)}$.

It can be seen that $\vec{r}_E$ is a residue system for $p$ by the same argument as the proof of Lemma \ref{lem:E''}.

The other properties follow straightforwardly.

\end{proof}

\begin{rmk}
    We may assume that if $\alpha\in\EE^p_\delta$, then there is a residue function $r:\P_\delta/p\to\P_\delta\cap M^\delta_\alpha$ that is order-preserving and satisfies that the residue $r(q)$ always extends the trace of $q$ to $M^\delta_\alpha$.
\end{rmk}

\begin{cor}
    For every $\delta<\kappa^+$, the poset $\P_\delta$ is strongly proper with respect to the models $(M^\delta_\alpha:\alpha\in\EE_\delta)$.
\end{cor}
\begin{proof}
    Let $p\in M^\delta_\alpha$. By Lemma \ref{lem:addingmodels} there is $q\leq p$ such that $\alpha\in\EE^q_\delta$. By Lemma \ref{lem:spconstres}, every extension $q'\leq q$ has a residue system $\vec{r}_E$ into $M^\delta_\alpha$. Its root condition $r_{(\delta,\alpha)}$ is a residue of $q'$ into $M^\delta_\alpha$.
\end{proof}

\begin{cor}
    The poset $\P_\delta$ is strongly proper with respect to stationarily many $M\in\Pow_\kappa(H_\lambda)$ for every large enough regular cardinal $\lambda$ and every $\delta<\kappa^+$.
\end{cor}
\begin{proof}
    There are stationarily many models $M$ in $\Pow_\kappa(H_\lambda)$ that are elementary in $H_\lambda$ and satisfy $M\cap V_\kappa=V_\alpha$ for some $\alpha\in\EE_\delta$. Fix one such model $M$. It follows from $M\cap V_\kappa=V_\alpha$ and the implicit identification of $\P_\delta$ with a subset of $V_\kappa$ that
    \[
    \P_\delta\cap M=\P_\delta\cap M^\delta_\alpha.
    \]
    Hence every strongly $(\P_\delta,M^\delta_\alpha)$-generic condition is $(\P_\delta,M)$-generic, and thus strong properness of $\P_\delta$ with respect to $M^\delta_\alpha$ implies strong properness of $\P_\delta$ with respect to $M$.
\end{proof}

\begin{cor} Let $G\subseteq\P_{\kappa^+}$ be a generic filter. 
    \begin{enumerate}
        \item $\aleph_1^{V[G]}=\kappa$,
        \item $\aleph_2^{V[G]}=\kappa^+$,
        \item $(2^\omega)^{V[G]}=(2^{\omega_1})^{V[G]}=\aleph_2^{V[G]}$.
    \end{enumerate}
\end{cor}
\begin{proof}
    The preservation of $\kappa$ follows from strong properness. Since every cardinal below $\kappa$ is collapsed to be countable, we must have $\aleph_1^{V[G]}=\kappa$. The preservation of $\kappa^+$ follows from $\kappa$-cc: the poset $\P_{\kappa^+}$ is a direct limit of posets that have $\kappa^+$-cc (for the trivial reasons that they have size $\kappa$), and thus $\P_{\kappa^+}$ must also have $\kappa^+$-cc. Hence $\aleph_2^{V[G]}=\kappa^+$. The third item follows also from strong properness; each $\P_\delta$, $\delta<\kappa^+$, has complete subposets of size $<\kappa$, and thus the continuum is pushed up to $\aleph_2^{V[G]}$. 
\end{proof}

\begin{cor}
    In $V^{\P_{\kappa^+}}$, the tree $\dot T$ is a wide $\kappa$-tree and for every $\gamma<\kappa^+$, there is an injective meet- and level-preserving tree-embedding
    \[
    \dot f_\gamma:\dot S_\gamma\to\dot T.
    \]
\end{cor}

It remains to show that the tree $\dot T$ remains Aronszajn throughout the iteration.

\section{Keeping $\dot T$ Aronszajn}

The last part of the proof is to show that the tree $\dot T$ is still Aronszajn in the final extension by $\P_{\kappa^+}$. Before the proof, we need to prove a quotient version of strong properness in order to carry out a splitting argument. By Lemma \ref{lem:spconstres} and its proof we assume that if $p\in\P_\delta$ and $\alpha\in\EE^p_\delta$, then there is an order-preserving residue map $\P_\delta/p\to\P_\delta\cap M^\delta_\alpha$.

\subsection{Residue systems in quotients}

\begin{notation}
    For $\delta\leq\bar{\delta}<\kappa^+$ and $p\in\P_\delta$, let 
    \[
    \EE^p_{\bar{\delta}}:=\{\alpha\in\EE_{\bar{\delta}}:\alpha\in\NN^p_\gamma\text{ for every }\gamma\in \delta\cap M^{\bar{\delta}}_{\alpha}\}.
    \]
\end{notation}

\begin{lem}\label{lem:pathcomp}
    Let $\delta\leq\bar{\delta}<\kappa^+$ and $p\in\P_\delta$. If $\bar{\alpha}\in\EE^{p}_{\bar{\delta}}$, then $p$ is strongly $(\P_\delta,M^{\bar{\delta}}_{\bar{\alpha}})$-generic.
\end{lem}
\begin{proof}
    We show that if $\bar{\alpha}\in\EE^{p}_{\bar{\delta}}$, then $p$ has a residue in $M^{\bar{\delta}}_{\bar{\alpha}}$. This implies the existence of a residue map $\P_\delta/p\to\P_\delta\cap M^{\bar{\delta}}_{\bar{\alpha}}$, since any extension $p'\leq p$ also must satisfy $\bar{\alpha}\in\EE^{p'}_{\bar{\delta}}$. To this end, note that $\bar{\alpha}\in\EE^{p}_{\bar{\delta}}$ implies the existence of a condition $q\in\P_{\bar{\delta}}$ such that $\bar{\alpha}\in\EE^q_{\bar{\delta}}$ and $q\rest\delta=p$. Indeed, $q$ can be defined by stretching $p$ by for every $\xi\in[\delta,\bar{\delta})$, letting $q(\xi)$ to be either $(\emptyset,\{\bar{\alpha}\})$ if $\xi\in M^{\bar{\delta}}_{\bar{\alpha}}$, and $q(\xi):=(\emptyset,\emptyset)$ if $\xi\notin M^{\bar{\delta}}_{\bar{\alpha}}$. Then $q\in\P_{\bar{\delta}}$ and $\bar{\alpha}\in\EE^q_{\bar{\delta}}$ so $q$ has a residue $r$ into $M^{\bar{\delta}}_{\bar{\alpha}}$, and $r\rest\delta$ is a residue for $p$ into $M^{\bar{\delta}}_{\bar{\alpha}}$: if $w\in \P_\delta\cap M^{\bar{\delta}}_{\bar{\alpha}}$ extends $r\rest\delta$, then $w^\smallfrown r\rest[\delta,\bar{\delta})$ extends $r$ in $\P_{\bar{\delta}}\cap M^{\bar{\delta}}_{\bar{\alpha}}$ and is compatible with $q$ and thus with $p$ too, since $r$ is a residue of $q$ into $M^{\bar{\delta}}_{\bar{\alpha}}$.
\end{proof}

If $G$ is a generic filter on a subposet of $\P$ and $p\in\P/G$, we say that a condition $r\in(\P/G)\cap M$ is a \textbf{residue of $\bs{p}$ into $\bs{M}$ in the quotient $\bs{\P/G}$} if for every $w\in(\P/G)\cap M$ that extends $r$ there is $q\in\P_\delta/G$ that extends $w$ and $p$. If there is an order-preserving residue map $\P/p\to\P\cap M$ and $G\subseteq\P\cap M$ is a generic filter that contains a residue of $p$, then $p\in\P/G$, and if $H$ is a $V[G]$-generic filter on $\P/G$, then in fact it is a $V$-generic filter on $\P$.

\begin{lem}
    If there is a path from $M^{\bar{\delta}}_{\bar{\alpha}}$ to $M^\delta_\alpha$, then $\P_\delta\cap M^{\bar{\delta}}_{\bar{\alpha}}\subseteq\P_\delta\cap M^\delta_\alpha$.
\end{lem}
\begin{proof}
    Follows from Lemma \ref{lem:claim52}.
\end{proof}

\begin{lem}\label{lem:resquo2}
    Let $\delta\leq\bar{\delta}<\kappa^+$, $\bar{\alpha}\leq\alpha<\kappa$ and $p\in\P_\delta$. Suppose that $\alpha\in\EE^p_\delta$, $\bar{\alpha}\in\EE^p_{\bar{\delta}}$ and there is a path from $M^{\bar{\delta}}_{\bar{\alpha}}$ to $M^\delta_\alpha$. Let $G$ be a generic filter on $\P_\delta\cap M^{\bar{\delta}}_{\bar{\alpha}}$ such that $p\in\P_\delta/G$. Assume that $\vec{r}_E$ is a residue system of $p$ into $M^\delta_\alpha$ whose root condition is in $\P_\delta/G$. Then for every $(\gamma,\beta)\in E$, the condition $r_{(\gamma,\beta)}$ is a residue for $r_{(\gamma,\beta^+)}$ into $M^\delta_\alpha$ in the quotient $\P_\delta/G$.
\end{lem}
\begin{proof}
    Note first that every member of the residue system must be in the quotient $\P_\delta/G$ because the root condition is in $\P/G$.
    Thus, arguing by induction, it suffices to show that the root condition $r_{(\delta,\alpha)}$ is a residue for $p$ into $M^\delta_\alpha$ in the quotient $\P_\delta/G$. Since $G$ contains a residue of $r_{(\delta,\alpha)}$ into $M^{\bar{\delta}}_{\bar{\alpha}}$, we have $r_{(\delta,\alpha)}\in(\P_\delta\cap M^\delta_\alpha)/G$. Let $w\in(\P_{\delta}/G)\cap M^\delta_\alpha$ extend $r_{(\delta,\alpha)}$. Then since $w$ is in $\P_{\delta}/G\cap M^\delta_\alpha=(\P_{\delta}\cap M^{\delta}_\alpha)/G$, there is a $V[G]$-generic filter $H$ on $(\P_{\delta}\cap M^{\delta}_\alpha)/G$ that contains $w$ and therefore $r_{(\delta,\alpha)}$ too. This $H$ is a $V$-generic filter on $\P_\delta\cap M^\delta_\alpha$ because there is an order-preserving residue map $(\P_\delta\cap M^\delta_\alpha)/r_{(\delta,\alpha)}\to\P_\delta\cap M^{\bar{\delta}}_{\bar{\alpha}}$. Also, it extends $G$. And since $H$ contains $r_{(\delta,\alpha)}$, which is a residue of $p$ into $M^{\delta}_\alpha$, we have $p\in\P_{\delta}/H$. So there is a generic filter $K$ on $\P_{\delta}/H$ that contains $p$. Now $K$ is a $V$-generic filter on $\P_{\delta}$ that contains $w$ and $p$. It follows that there is a common extension of $w$ and $p$ in $H$, and this common extension must be in $\P_{\delta}/G$ because $H$ extends $G$.
\end{proof}

\begin{lem}\label{lem:spquo}
    Let $\delta\leq\bar{\delta}<\kappa^+$, $\bar{\alpha}\leq\alpha<\kappa$ and $p\in\P_\delta$. Suppose that $\alpha\in\EE^p_\delta$, $\bar{\alpha}\in\EE^p_{\bar{\delta}}$ and there is a path from $M^{\bar{\delta}}_{\bar{\alpha}}$ to $M^\delta_\alpha$. Let $G$ be a generic filter on $\P_\delta\cap M^{\bar{\delta}}_{\bar{\alpha}}$ such that $p\in\P_\delta/G$. Then $p$ has a residue system into $M^\delta_\alpha$ whose every condition is in $\P_\delta/G$.
\end{lem}

\begin{proof}
    The proof is by induction on $\delta$. We look at the successor case $\delta+1$.
    
    Let $H\subseteq\P_{\delta+1}/G$ be a $V[G]$-generic filter that contains $p$. Then $H$ is a $V$-generic filter on $\P_\delta$ that extends $G$. For any $\gamma\leq\delta+1$ and $\beta\leq \kappa$, denote
    \[
    H^\gamma_\beta:=H\cap\P_\gamma\cap M^\gamma_\beta.
    \]
    Enumerate the set $\EE^p_\delta-\alpha$ as
    \[
    \alpha=\beta_0<\dots<\beta_{n-1}.
    \]
    Abbreviate $\beta_{n}:=\kappa$ and $r_{n}:=p\rest\delta$. Now
    \[
    H^\delta_\alpha=H^\delta_{\beta_0}\subseteq\dots\subseteq H^\delta_{\beta_{n+1}}
    \]
    and each $H^\delta_{\beta_k}$ is a $V$-generic filter on $\P_\delta\cap M^\delta_{\beta_k}$.
    We proceed by reverse recursion on $k<n$.
    Assume that we have defined $r_{k+1}$ and $r_{k+1}\in \P_\delta/H^\delta_{\beta_{k}}$.
    By the induction hypothesis $r_{k+1}$ has a residue system $\vec{r}_{E_k}$ into $M^\delta_{\beta_k}$ with every condition in $\P_\delta/H^\delta_{\beta_k}$. In particular, the root condition $r_k:=r_{(\delta,\beta_k)}$ is an element in $H^\delta_{\beta_k}$. 
    Let $X_k$ be the collection of pairs $(s,t)\in f^p_\delta$ of exit nodes from $V_\alpha$ such that $t\in M^\delta_{\beta_{k+1}}- M^\delta_{\beta_k}$. For each pair $(s,t)\in X_k$, let $t_{\bar{s}}:=\pi^{E_k}(t)$.
    Then $t_{\bar{s}}$ is a node in $T\cap V_{\alpha}$ since $t$ is an exit node from $V_\alpha$.
    The branch below $s$ is decided by the poset $\P_\delta\cap M^\delta_{\beta_k}$, because $\P_\delta,\dot S_\delta,s\in M^\delta_{\beta_k}$. Therefore the branch below $s$ is in $V[H^\delta_{\beta_k}]$. Let $\bar{s}$ be the predecessor of $s$ at the height of the node $t_{\bar{s}}$, as decided by the generic $H^\delta_{\beta_k}$. Up to extending $r_k$ inside the generic $H^\delta_{\beta_k}$, we may assume that it forces that $\bar{s}$ is the implicit preimage of $t_{\bar{s}}$, for each such pair $(s,t)$. Then, since $r_k\in H^\delta_{\beta_k}$, in particular also $r_k\in\P_\delta/H^\delta_{\beta_{k-1}}$ and $r_k\in\P_\delta/G$.

    Finally, suppose that we have defined $r_1$. Then $r_1\in \P_\delta/H^\delta_{\beta_0}$. By the induction hypothesis $r_1$ has a residue system $\vec{r}_{E_1}$ into $M^\delta_\alpha$ with root condition $r_0$ in $\P_\delta/H^\delta_\alpha$. Define
    \[
    f:=(f^p_\delta\cap V_\alpha)\cup\bigcup_{k\leq n}\{(\bar{s},t_{\bar{s}}):(s,t)\in X_k\},
    \]
    let $\NN:=\NN^p_\delta\cap\alpha$ and let
    \[
    r:={r_0}^\smallfrown(f,\NN).
    \]
    We claim that $r\in H$. Note that $r_0\in H$. Note also that $f$ is a finite level- and meet-preserving function from $\dot S_\delta^{H^\delta_\alpha}$ to $\dot T^{H^\delta_\alpha}$. Thus, it suffices to show that $f\subseteq f^H_\delta$ and $\NN\subseteq\NN^G_\delta$, where
    \begin{align*}
        &f^H_\delta=\bigcup_{q\in H}f^q_\delta,\\
        &\NN^H_\delta=\bigcup_{q\in H}\NN^q_\delta.
    \end{align*}
    The fact that $\NN\subseteq\NN^H_\delta$ follows immediately from the fact that $p\in H$ and $\NN=\NN^p_\delta\cap\alpha$. We show that $f\subseteq f^H_\delta$.
    Note first that the function $f^H_\delta$ is an injective level- and meet-preserving tree-embedding from $\dot S^H_\delta$ into $\dot T^H_\delta$. Furthermore, since $p\in H$, it extends the finite partial function $f^p_\delta$.
    Furthermore, if $(s,t)\in f^H_\delta$ and $\bar{s}<_{\dot S^H_\delta}s$ and $\bar{t}<_{\dot T^H}t$ are such that $\height(\bar{s})=\height(\bar{t})$, then also $(\bar{s},\bar{t})\in f^H_\delta$. Now, if $(\bar{s},\bar{t})\in f$, then there are $(s,t)\in f^p_\delta\subseteq f^H_\delta$ such that $r_k\Vdash ``\bar{s}<s$ and $\bar{t}<t"$. Since $r_k\in H$ for each $k$, it follows that $f\subseteq f^H_\delta$.

    We have thus shown that $r\in H$. Since $r\in M^{\delta+1}$, it follows that $r\in(\P_{\delta+1}\cap M^{\delta+1}_\alpha)/G$. The system $\vec{r}_E$ where $E=\bigcup_{k\leq n}\cup\{(\delta+1,\alpha)\}$ and $r_{(\delta+1,\alpha)}:=r$ is thus a residue system for $p$ into $M^{\delta+1}_\alpha$ with every condition in the quotient $\P_{\delta+1}/G$.
\end{proof}

\vv

\subsection{No new branches}\spa

\vv

What follows in this section is heavily inspired by the argument in Laver-Shelah \cite{laver1981aleph} that allows one to amalgamate partial specialization functions using ``splitting".

\vv

For a model $M\elem H_\theta$ and a generic filter $G$ on a poset $\P\in M$, we denote $M[G]=\{\dot a^G:\dot a\in V^{\P}\cap M\}$.

\begin{lem}\label{lem:TV}
    If $G$ is a generic filter on $\P_\delta$ and $\alpha\in\EE_\delta$, then $M^\delta_\alpha[G]$ is closed under $<\alpha$-sequences and elementary in $H_\theta[G]$.
\end{lem}
\begin{proof}
    The fact that $M^\delta_\alpha[G]$ is elementary in $H_\theta[G]$ follows using the Tarski-Vaught criterion and the Maximality principle: if $H_\theta[G]\models\exists x\phi(x,\dot a^G)$ for some $\dot a^G\in M[G]$, then there is $\dot b$ such that $\Vdash\exists x\phi(x,\dot a)\to\phi(\dot b,\dot a)$, and by elementarity $M^\delta_\alpha\elem H_\theta$, such $\dot b$ can be found in $M^\delta_\alpha$. Then $H_\theta[G]\models\phi(\dot b^G,\dot s^G)$ where $\dot b^G\in M^\delta_\alpha[G]$ and we are done.

    To see that $M^\delta_\alpha[G]$ is closed under $<\alpha$-sequences, let $\bar{\alpha}<\alpha$ and let $a:\bar{\alpha}\to M^\delta_\alpha[G]$ be a function. We show that $a\in M^\delta_\alpha[G]$. For each $i<\bar{\alpha}$, let $\dot a_i\in V^{\P_\delta}\cap M^\delta_\alpha$ be a name for $a(i)$. Since $M^\delta_\alpha$ is closed for $<\alpha$-sequences, the $\P_\delta$-name
    \[
    \{((i,\dot a_i),1_{\P_\delta}):i<\bar{\alpha}\}
    \]
    is a name for $a$ in $M^\delta_\alpha$ and so $a=\dot a^G\in M^\delta_\alpha[G]$.
\end{proof}

The weakly compact cardinal is used in the following lemma. We suppose without loss of generality that the name $\dot S$ in the statement of the lemma is a $\P_\delta$-name for a normal wide $\kappa$-Aronszajn tree on $\kappa\times\kappa$ such that $\Vdash_{\P_\delta}`` \Lev_\beta(\dot S)=\kappa\times\{\beta\}"$ for every $\beta<\kappa$.

In the following definition, the case $s=s'$ is not excluded.

\begin{de}
    Let $\P$ be a poset and let $\dot S$ be a $\P$-name for a tree. A pair of conditions $(p,q)$ \textbf{splits a pair of nodes $\bs{(s,s')}$} if there are distinct nodes $\bar{s}$ and $\bar{s}'$ and an ordinal $\bar{\alpha}$ such that
    \begin{enumerate}
        \item $p\Vdash``\bar{s}<s$ and $\height(\bar{s})=\bar{\alpha}"$,
        \item $q\Vdash``\bar{s}'<s'$ and $\height(\bar{s}')=\bar{\alpha}"$.
    \end{enumerate}
    In this case we say that $p$ and $q$ split the pair $(s,s')$ with the pair $(\bar{s},\bar{s}')$.
\end{de}

\begin{lem}\label{lem:split:onestep} Let $\delta<\kappa^+$ and $\alpha\in\EE_\delta$. Suppose that $\dot S\in M^\delta_\alpha$ is a $\P_\delta$-name for a wide $\kappa$-Aronszajn tree and $G\subseteq\P_\delta\cap M^\delta_\alpha$ is a generic filter. Suppose that $p$ and $q$ are two conditions in the quotient $\P_\delta/G$. Then for any nodes $s,s'\in \dot S$ of limit height that are an exit nodes from $V_\alpha$ there are two extensions $\hat{p}\leq p$ and $\hat{q}\leq q$ in the quotient $\P_\delta/G$ that split the pair $(s,s')$.
\end{lem}

\begin{proof}
    Let $s,s'\in\dot S$ be nodes at a limit level that are exit nodes from $V_\alpha$. We first claim that there are two conditions $p^L,p^R\leq p$ in $\P_\delta/G$ that split $s$. If not, then there is $b\in V[G]$ such that $p\Vdash ``b$ is the branch below $s$". Let $\bar{\alpha}$ be the height of $s$ and let $\bar{\beta}$ be the width of $s$, i.e. the least $\beta\leq\kappa$ such that $b\subseteq\bar{\alpha}\times\beta$. We consider three cases.

    \begin{itemize}
    \item \textbf{Case $\bs{1}$:} $\bar{\alpha}=\alpha$.
    Now the branch $b$ is a cofinal branch in the tree $(\dot S\cap V_\alpha)^G$. But since $M^\delta_\alpha$ reflects all $\Pi^1_1$-statements with parameters in it, and indeed $\dot S,\P_\delta\in M^\delta_\alpha$, the tree $(\dot S\cap V_\alpha)^G$ must be a wide $\alpha$-Aronszajn tree. Hence it is impossible to have $b$ in $V[G]$.
        
    \item \textbf{Case $\bs{2}$:} $\bar{\alpha}<\alpha$ and $\bar{\beta}=\alpha$.
    Now the branch $b$ induces a cofinal function from $\bar{\alpha}$ to $\bar{\beta}$ in $V[G]$, by the definition of $\bar{\beta}$. But this is not possible, because $\bar{\beta}=\alpha=\omega_2^{V[G]}$ and $\bar{\alpha}<\bar{\beta}$.
    
    \item \textbf{Case $\bs{3}$:} $\bar{\alpha}<\alpha$ and $\bar{\beta}<\alpha$.
    Let $G^+$ be a generic filter on $\P_\delta$ that extends $G$ and contains $p$. Now $b$ is a bounded subset of $M^\delta_\alpha[G^+]$ and so $b\in M^\delta_\alpha[G^+]$, by Lemma \ref{lem:TV}. Furthermore $H_\theta[G^+]$ satisfies that the branch $b$ has a supremum in $\dot S^{G^+}$, namely the node $s$. Since $M^\delta_\alpha[G^+]\elem H_\theta[G^+]$ (again by Lemma \ref{lem:TV}), the model $M^\delta_\alpha[G^+]$ must also contain a supremum of $b$. But this is impossible since the unique supremum of the branch $b$ is the node $s$, and $s$ is an exit node from $M^\delta_\alpha$, and thus also from $M^\delta_\alpha[G^+]$.
    \end{itemize}
Hence we may find two extensions $p^L,p^R\leq p$ in the quotient $\P_\delta/G$ that split $s$ at some height $\bar{\alpha}<\alpha$ with some distinct nodes $s^L$ and $s^R$. In the quotient $\P_\delta/G$, find an extension $\hat{q}\leq q$ that decides the predecessor of $s'$ at height $\bar{\alpha}$, call it $\bar{s}'$. If $s^L\neq\bar{s}'$, let $\hat{p}:=p^L$, and otherwise let $\hat{p}:=p^R$. Then $\hat{p}$ and $\hat{q}$ are as wanted.

\end{proof}

\begin{lem}\label{lem:succ}
    Let $\delta<\kappa^+$, $\alpha\in \EE_\delta$ and let $G\subseteq\P_\delta\cap M^\delta_\alpha$ be a generic filter. Assume that $r^L,r^R\in\P_\delta/G$ are two conditions and $f$ and $N$ satisfy
    \begin{enumerate}
        \item ${r^L}^\smallfrown(f,N)$ and ${r^R}^\smallfrown(f,N)$ are conditions in $\P_{\delta+1}$ and $\alpha\in N$,
        \item $r^L$ and $r^R$ decide the widths of nodes in the set $\dom(f)$ similarly,
        \item the conditions $r^L$ and $r^R$ split any pair of nodes from $\dom(f)$ that are exit nodes from $V_\alpha$ at a limit level with some nodes in $\dom(f)\cap V_\alpha$,
        \item if $s\in\dom(f)$ is an exit node from $V_\alpha$ at a successor level, then $r^L$ and $r^R$ decide the immediate predecessor of $s$ the same way and it is in $\dom(f)\cap V_\alpha$.
    \end{enumerate}
    Then the conditions ${r^L}^\smallfrown(f,N)$ and ${r^R}^\smallfrown(f,N)$ have a common residue into $M^{\delta+1}_\alpha$.
\end{lem}
\begin{proof}
    Note first that the function $f\cap V_\alpha$ is a level- and meet-preserving tree-embedding from $(\dot S_\delta\cap V_\alpha)^G$ to $(\dot T\cap V_\alpha)^G$, whose domain is closed under meets.
    
    We build residue systems for ${r^L}^\smallfrown(f,N)$ and ${r^R}^\smallfrown(f,N)$ by reverse recursion on $k\leq n$. In the beginning, let $v^L_n:=r^L$ and $v^R_n:=r^R$. Then $v^L_n,v^R_n\in\P_\delta/G$.

    At step $k<n$, suppose that $v^L_{k+1}$ and $v^R_{k+1}$ are defined and are in $(\P_\delta\cap M^\delta_{\beta_{k+1}})/G$. By Lemma \ref{lem:spquo}, they have residue systems $\vec{v}^L_{F_k^L}$ and $\vec{v}^R_{F^R_k}$ into $M^\delta_{\beta_k}$ whose root conditions $v^L_k$ and $v^R_k$ are in the quotient $\P_\delta/G$. These residue systems can be chosen inside the model $M^\delta_{\beta_{k+1}}$. Exactly as in the proof of Lemma \ref{lem:spconstres} we extend the root conditions $v^L_k$ and $v^R_k$ inside the model $M^\delta_{\beta_k}$ to decide, for every node $s$ in the set
    \[
    X_k:=\{s\in\dom(f):s\text{ is an exit node from }V_\alpha\text{ and }f(s)\in M^\delta_{\beta_{k+1}}-M^\delta_{\beta_k}\},
    \]
    the implicit preimages $\bar{s}^L$ and $\bar{s}^R$ for the nodes $t_{\bar{s}^L}:=\pi^{F^L_k}(f(s))$ and $t_{\bar{s}^R}:=\pi^{F_k^R}(f(s))$, respectively. Note that for any two nodes $s,s'\in X_k$ (possibly $s=s'$), if $\bar{s},\bar{s}'\in\dom(f)$ are distinct nodes such that
    \[
    r^L\Vdash \bar{s}^L<_{\dot S_\delta}s\quad\text{and}\quad r^R\Vdash\bar{s};<_{\dot S_\delta}s',
    \]
    it then holds that 
    \[
    \bar{s}\leq_{(\dot S_\delta\cap V_\alpha)^G}s\quad\text{and}\quad \bar{s}\leq_{(\dot S_\delta\cap V_\alpha)^G}s'.
    \]
    In particular the meet of $\bar{s}^L$ and $\bar{s}^R$ is the meet of $\bar{s}$ and $\bar{s}'$, and thus already an element in $\dom(f)$. Furthermore, the meet of $t_{\bar{s}^L}$ and $t_{\bar{s}'^R}$ must be the meet of the nodes $f(\bar{s})$ and $f(\bar{s}')$. So in particular the function
    \[
    f_k:=(f\cap V_\alpha)\cup\{(\bar{s}^L,t_{\bar{s}^L}),(\bar{s}^R,t_{\bar{s}^R}):s\in\bigcup_{k\leq l<n}X_l\}
    \]
    must be a level- and meet-preserving tree-embedding from $(\dot S_\delta\cap V_\alpha)^G$ to $(\dot T\cap V_\alpha)^G$ whose domain is closed under meets.

    Finally, after the recursion, the function $f_0$ is a finite level- and meet-preserving tree-embedding from $(\dot S\cap V_\alpha)^G$ to $(\dot T\cap V_\alpha)^G$ that extends $f\cap V_\alpha$. Without loss of generality, $v_0:=v^L_0=v^R_0\in G$, and the concatenation ${v_0}^\smallfrown(\hat{f},\NN^p_\delta\cap\alpha)$ is a condition in $(\P_{\delta+1}\cap M^{\delta+1}_\alpha)/G$. Let $F^L$ be the union of the sets $F_k^L$, $k\leq n$, together with the pair $(\delta+1,\alpha)$ and similarly let $F^R$ be the union of the $F^R_k$, $k\leq n$, together with the pair $(\delta+1,\alpha)$. Let
    \begin{align*}
        v_{(\delta+1,\alpha)}^L=v^R_{(\delta+1,\alpha)}={v_0}^\smallfrown(f_0,N\cap\alpha)
    \end{align*}
    Then $\vec{v}^L_{F^L}$ and $\vec{v}^R_{F^R}$
    are residue systems for ${r^L}^\smallfrown(f,N)$ and ${r^R}^\smallfrown(f,N)$, respectively, into $M^{\delta+1}_\alpha$, with the common root condition. The common root condition is a common residue for ${r^L}^\smallfrown(f,N)$ and ${r^R}^\smallfrown(f,N)$ into $M^{\delta+1}_\alpha$.
\end{proof}

\begin{lem}\label{lem:splitbiglem}
    Assume that the bookkeeping function is such that $\dot S_\gamma$ is a $\P_\gamma$-name for a wide $\kappa$-Aronszajn tree, for every $\gamma<\kappa^+$. Let $\delta<\kappa^+$. 
    \begin{enumerate}
        \item\label{item1:lem:splitbiglem} For any $p\in\P_\delta$, $\alpha\in\EE^p_\delta$ and any finite set of nodes $A\subseteq\dot T$ of limit height that are exit nodes from $M^\delta_\alpha$, there are two conditions $q^L,q^R\leq p$ that have a common residue into $M^\delta_\alpha$ and split every node in $A$.
        \item\label{item2:lem:splitbiglem} $\dot T$ is a $\kappa$-Aronszajn tree in $V^{\P_\delta}$.
    \end{enumerate}
\end{lem}
\begin{proof}
    The proof is by induction on $\delta$. We prove item (\ref{item1:lem:splitbiglem}) first. The base case follows using Proposition \ref{prop:1tree}(\ref{prop:1tree:free}). The limit case follows straightforwardly from the induction hypothesis. We consider the successor case $\delta+1$.
    
    Let $p\in\P_{\delta+1}$, $\alpha\in\EE^p_{\delta+1}$ and let $A\subseteq\dot T$ be a finite set of exit nodes from $M^{\delta+1}_\alpha$. Suppose up to extending $p$ that whenever $s\in\dom(f^p_\delta)$ is an exit node from $V_\alpha$ at a successor height, then $p\rest\delta$ decides its immediate predecessor and this predecessor is in $\dom(f^p_\delta)$. Let $\vec{r}_E$ be a residue system for $p$ into $M^{\delta+1}_\alpha$. Enumerate the set $E_\delta$ as $\beta_0<\dots<\beta_{n-1}$, abbreviate $r_k:=r_{(\delta,\beta_k)}$, denote $r_{n}:=p\rest\delta$, $\beta_{n}:=\kappa$ and let
    \[
    E_k:=E\cap(M^\delta_{\beta_{k+1}}-M^\delta_{\beta_k}).
    \]
    Then $\vec{r}_{E_k}$ is a residue system for $r_{k+1}$ into $M^\delta_{\beta_k}$, for every $k< n$, by Lemma \ref{lem:E''}. Let $G\subseteq\P_\delta\cap M^\delta_\alpha$ be a generic filter that contains the condition $r_0$. Then $r_k\in(\P_\delta\cap M^\delta_{\beta_k})/G$ for each $k$ and $r_{(\delta+1,\alpha)}\in (\P_{\delta+1}\cap M^{\delta+1}_\alpha)/G$. For simplicity, we may suppose that $f^{r_{(\delta+1,\alpha)}}_\delta$ and $\NN^{r_{(\delta+1,\alpha)}}_\delta$ are as in the proof of Lemma \ref{lem:spconstres}; $f^{r_{(\delta+1,\alpha)}}_\delta$ is obtained from $f^p_\delta$ by adding pairs $(\bar{s},t_{\bar{s}})$ obtained from pairs $(s,t)\in f^p_\delta$ of exit nodes from $V_\alpha$ by letting $t_{\bar{s}}:=\pi^{E_k}(t)$, where $k$ is such that $t\in M^\delta_{\beta_{k+1}}-M^\delta_{\beta_k}$, and $\bar{s}$ is forced by $r_k$ to be the predecessor of $s$ at the height of $t_{\bar{s}}$, and $\NN^{r_{(\delta+1,\alpha)}}_\delta=\NN^p_\delta\cap\alpha$. It holds in particular that if $r$ is any common extension of $r_0,\dots,r_n$ and $p\rest\delta$ in $\P_\delta$, then $r^\smallfrown(f^{r_{(\delta+1,\alpha)}}_\delta\cup f^p_\delta,\NN^p_\delta)$ is a condition in $\P_{\delta+1}$. For each $k<n$, let
    \[
    X_{k}:=\{s\in\dom(f^p_\delta):s\text{ is an exit node from }V_\alpha\text{ and }f^p_\delta(s)\in(M^\delta_{\beta_{k+1}}-M^\delta_{\beta_k})\}.
    \]
    Again, $X_k\subseteq M^\delta_{\beta_k}$ and if $s\in X_k$, then $\lab(f^p_\delta(s))\leq\beta_k$. We build conditions $r^L_k,r^R_k\leq r_k$ in $(\P_\delta\cap M^\delta_{\beta_k})/G$ and functions $f_k$ by recursion on $k\leq n+1$. 
    
    Let $r^L_0=r^R_0:=r_0$ and $f_0:=f^{r_{(\delta+1,\alpha)}}_\delta$. Suppose that $r^L_k$, $r^R_k$ and $f_k$ were defined. Suppose that they satisfy:
    \begin{enumerate}
        \item $r^L_k,r^R_k\in(\P_\delta\cap M^\delta_{\beta_k})/G$ and $r^L_k,r^R_k\leq r_k$,
        \item $r^L_k$ and $r^R_k$ split every node in $A\cap M^\delta_{\beta_k}$,
        \item $r^L_k$ and $r^R_k$ split every pair of nodes $(s,s')$ from $\bigcup_{j\leq k}X_j$, with a pair nodes $(\bar{s},\bar{s}')$,
        \item $f_k$ is a level- and meet-preserving tree-embedding from $(\dot S_\delta\cap V_\alpha)^G$ to $(\dot T\cap V_\alpha)^G$ such that extends $f^{r_{(\delta+1,\alpha)}}_\delta$ and the domain satisfies
        \[
        \dom(f^{r_{(\delta+1,\alpha)}}_\delta)\cup\{\bar{s},\bar{s}':(s,s')\in(\bigcup_{j\leq k}X_j)^2\}.
        \]
        and both $r^L_k$ and $r^R_k$ force that $f_k\cup f^p_\delta$ is a level- and meet-preserving tree-embedding from $\dot S_\delta$ to $\dot T$.
    \end{enumerate}
    We build $r^L_{k+1}$, $r^R_{k+1}$ and $f_{k+1}$.

    \begin{claim}
        There are common extensions $r^L_{k+1}\leq r^L_k,r_{k+1}$ and $r^R_{k+1}\leq r^R_k,r_{k+1}$ in the quotient $(\P_\delta\cap M^\delta_{\beta_{k+1}})/G$ such that for every pair $(s,s')\in \bigcup_{j\leq k}X_j$,
        \begin{align*}
            &r^L_{k+1}\Vdash f_k(\bar{s})<_{\dot T} f^p_\delta(s),\\
            &r^R_{k+1}\Vdash f_k(\bar{s}')<_{\dot T}f^p_\delta(s').
        \end{align*}
    \end{claim}
    \begin{proof}
        Using Lemma \ref{lem:freemodifs}, as in the proof of Lemma \ref{lem:spsysisres}. We get the common extensions in the quotient $(\P_\delta\cap M^\delta_{\beta_{k+1}})/G$ by Lemma \ref{lem:resquo2}.
    \end{proof}

    \noindent It now holds that both $r^L_{k+1}$ and $r^R_{k+1}$ force that the function
    \[
    f_k\cup(f^p_\delta\cap M^\delta_{\beta_{k+1}})
    \]
    is a level- and meet-preserving tree-embedding.

    \begin{claim}\label{claim:split:extended}
        There are extensions $\tilde r^L_{k+1}\leq r^L_{k+1}$ and $\tilde r^R_{k+1}\leq r^R_{k+1}$ in $(\P_\delta\cap M^\delta_{\beta_{k+1}})/G$ that split:
        \begin{enumerate}
            \item\label{item:type1} every node in $A\cap M^\delta_{\beta_{k+1}}$, and
            \item\label{item:type2} every pair of nodes $(s,s')$ from $\bigcup_{j\leq k+1}X_{j}$ with distinct nodes $(\bar{s},\bar{s}')$.
        \end{enumerate}
    \end{claim}
    \begin{proof}
        Using Lemma \ref{lem:split:onestep} repeatedly.
    \end{proof}

    Up to extending $r^L_k$ and $r^R_k$ in the quotient, we assume that they satisfy Claim \ref{claim:split:extended}. So any pair $(s,s')\in \bigcup_{j\leq k}X_j$ is split by a pair $(\bar{s},\bar{s}')$ by $r^L_k$ and $r^R_k$. Note that the meet $\bar{s}\wedge \bar{s}'$ is calculated in $V[G]$ in the tree $(\dot S_\delta\cap V_\alpha)^G$. It holds that:
    \begin{itemize}
        \item $r^L_{k+1}$ forces $s\wedge \bar{s}'=\bar{s}\wedge \bar{s}'$, and
        \item $r^R_{k+1}$ forces $\bar{s}\wedge s'=\bar{s}\wedge \bar{s}'$.
    \end{itemize}
    It follows that both $r^L_{k+1}$ and $r^R_{k+1}$ force that the set 
    \[
    (\dom(f^p_\delta)\cap V_{\beta_{k+1}})\cup\{\bar{s},\bar{s}':(s,s')\in (\bigcup_{j\leq k}X_j)^2\}
    \]
    is closed under meets. Consider the function $f_k$. By assumption it is a level- and meet-preserving tree-embedding from $(\dot S_\delta\cap V_\alpha)^G$ to $(\dot T\cap V_\alpha)^G$. By the Node Density Lemma \ref{lem:dens} applied in $(\P_{\delta+1}\cap M^{\delta+1}_\alpha)/G$, there is a level- and meet-preserving tree-embedding $f_{k+1}\supseteq f_k$ whose domain is the set
    \[
    \dom(f_k)\cup\{\bar{s},\bar{s}':(s,s')\in(\bigcup_{j\leq k}X_k)^2\},
    \]
    Indeed, the domain of $f_{k+1}$ must be closed under meets because the meet of $\bar{s}$ and $\bar{s}'$ in $(\dot S_\delta\cap V_\alpha)^G$ is the meet of $s$ and $s'$, as $p\rest\delta$, $r^L_k$ and $r^R_k$ all decide it, which must already be in the domain of the function $f^p_\delta\cap V_\alpha$. This ends the definition of $r^L_{k+1}$, $r^R_{k+1}$ and $f_{k+1}$.

    \vv

    Finally, suppose that we have defined $r^L_n$ and $r^R_n$. They are in $\P_\delta/G$ and extend each $r_k$, $k\leq n$. In particular, they extend $p\rest\delta$. Define
    \[
    f:=f^p_\delta\cup f_n.
    \]
    By construction, both $r^L_n$ and $r^R_n$ force that $f$ is a level- and meet-preserving tree-embedding from $\dot S_\delta$ to $\dot T$. Let
    \begin{align*}
        &q^L:={r^L_n}^\smallfrown(f,\NN^p_\delta),\\
        &q^R:={r^R_n}^\smallfrown(f,\NN^p_\delta).
    \end{align*}
    Furthermore, by construction, both $q^L$ and $q^R$ are conditions in $\P_{\delta+1}$. They split every node in $A$. They also split every pair $(s,s')$ of exit nodes from $V_\alpha$ in $f$. By Lemma \ref{lem:succ} $q^L$ and $q^R$ have residue systems $\vec{v}^L_{F^L}$ and $\vec{v}^R_{F^R}$ into $M^{\delta+1}_\alpha$ that have the same root condition $v^L_{(\delta+1,\alpha)}=v^R_{(\delta+1,\alpha)}$.

    We prove item (\ref{item2:lem:splitbiglem}). Let $\delta<\kappa^+$ and suppose that the lemma holds for every $\gamma<\delta$. By what we just proved, item (\ref{item1:lem:splitbiglem}) holds for $\delta$. Suppose to the contrary that there are $\dot b\in V^{\P_{\delta}}$ and $p\in\P_{\delta}$ such that
    \[
    p\Vdash ``\dot b\text{ is a cofinal branch in }\dot T".
    \]
    Let $\alpha\in\EE_{\delta}$ be such that $p,\dot b\in M^{\delta}_\alpha$. Let $q\leq p$ be such that $\alpha\in\EE^q_{\delta}$. Up to extending $q$, assume that it decides the node $t:=\dot b(\alpha)$. By item (\ref{item1:lem:splitbiglem}) there are two conditions $q^L,q^R\leq q$ that split $t$ with some distinct nodes $t^L$ and $t^R$ at some height $\bar{\alpha}<\alpha$, and have a common residue into $M^{\delta}_\alpha$. Let $r$ be the common residue. Find $w\in\P_{\delta+1}\cap M^{\delta+1}_\alpha$ that extends $r$ and decides the node $\bar{t}:=\dot b(\bar{\alpha})$. Then, if $\bar{t}\neq t^L$, $w$ cannot be compatible with $q^L$, and if $\bar{t}\neq t^R$, $w$ cannot be compatible with $q^R$. This contradicts the fact that $r$ is a common residue for $q^L$ and $q^R$, and ends the proof of item (\ref{item2:lem:splitbiglem}).

\end{proof}

\section{Conclusion}

\noindent We may conclude:

\begin{thm} Assume that there exists a weakly compact cardinal. It is consistent that there is a universal wide $\aleph_1$-Aronszajn tree.
\end{thm}
\begin{proof}
    Let $\kappa\to H(\kappa^+)$, $\gamma\mapsto\dot S_\gamma$ be a bookkeeping function such that for every $\P_{\kappa^+}$-name $\dot S$ for a wide $\kappa$-Aronszajn tree there is $\gamma<\kappa^+$ such that $\dot S_\gamma$ is a $\P_\gamma$-name for a wide $\kappa$-Aronszajn tree and $\Vdash_{\P_\gamma}`` \dot S_\gamma\cong\dot S"$. Such a function exists by the fact that $\P_{\kappa^+}$ has $\kappa^+$-cc. Then, in $V^{\P_{\kappa^+}}$, $\aleph_1$ is $\kappa$, and $\dot T$ is a wide $\aleph_1$-tree that is universal for all wide $\aleph_1$-Aronszajn trees. Also $\dot T$ must be Aronszajn, since otherwise by $\kappa^+$-cc there is $\delta<\kappa^+$ such that $\dot T$ is not Aronszajn in $V^{\P_\delta}$, which contradicts Lemma \ref{lem:splitbiglem}. 
\end{proof}

I conjecture that the forcing construction can be run over $L$ using the combinatorial principle $\Diamond^{\#}$ from \cite{devlin1982combinatorial} and the filter derived from it, instead of collapsing a weakly compact cardinal and using the weakly compact filter. This would mean that the existence of a universal wide $\aleph_1$-Aronszajn tree could be shown to be consistent without any large cardinal assumptions.

\bibliographystyle{plain}
\bibliography{bib.bib}

\begin{thebibliography}{10}

\bibitem{ben2023aronszajn}
Omer Ben-Neria, Siiri Kivim{\"a}ki, Menachem Magidor, and Jouko V{\"a}{\"a}n{\"a}nen.
\newblock Aronszajn trees and maximality.
\newblock {\em \href{https://arxiv.org/abs/2305.07880}{arXiv:2305.07880}}, 2024.
\newblock Submitted.

\bibitem{devlin1982combinatorial}
Keith~J Devlin.
\newblock The combinatorial principle {$\Diamond^{\#}$}.
\newblock {\em The Journal of Symbolic Logic}, 47(4):888--899, 1982.

\bibitem{dvzamonja2006properties}
Mirna D{\v{z}}amonja and Saharon Shelah.
\newblock On properties of theories which preclude the existence of universal models.
\newblock {\em Annals of Pure and Applied Logic}, 139(1-3):280--302, 2006.

\bibitem{dvzamonja2021wide}
Mirna D{\v{z}}amonja and Saharon Shelah.
\newblock On wide {A}ronszajn trees in the presence of {MA}.
\newblock {\em The Journal of Symbolic Logic}, 86(1):210--223, 2021.

\bibitem{gilton2017side}
Thomas Gilton and Itay Neeman.
\newblock Side conditions and iteration theorems.
\newblock {\em Appalachian Set Theory}, 2017.

\bibitem{hyttinen1987games}
Tapani Hyttinen.
\newblock {\em Games and infinitary languages}.
\newblock Ann. Acad. Sci. Fenn. Ser. A I Math. Dissertationes No 64, 1987.

\bibitem{hyttinen1990scott}
Tapani Hyttinen and Jouko V{\"a}{\"a}n{\"a}nen.
\newblock On {S}cott and {K}arp trees of uncountable models.
\newblock {\em The Journal of symbolic logic}, 55(3):897--908, 1990.

\bibitem{jech2003set}
Thomas Jech.
\newblock {\em Set theory: The third millennium edition, revised and expanded}.
\newblock Springer, 2003.

\bibitem{kanamori2008higher}
Akihiro Kanamori.
\newblock {\em The higher infinite: large cardinals in set theory from their beginnings}.
\newblock Springer Science \& Business Media, 2008.

\bibitem{Karp1965-KARFQE}
Carol Karp.
\newblock Finite quantifier equivalence.
\newblock In J.~W. Addison, editor, {\em Journal of Symbolic Logic}, pages 407--412. Amsterdam: North-Holland Pub. Co., 1965.

\bibitem{karttunen1984model}
Maaret Karttunen.
\newblock {\em Model theory for infinitely deep languages}.
\newblock Add. Acad. Sci. Fenn. Ser. A I Math. Dissertationes No. 50, 1984.

\bibitem{kurepa1956ensembles}
Georges Kurepa.
\newblock Ensembles ordonn{\'e}s et leurs sous-ensembles bien ordonn{\'e}s.
\newblock {\em Comptes rendus hebdomadaires des s{\' e}ances de l'{A}cad{\'e}mie des sciences}, 242(18):2202--2203, 1956.

\bibitem{laver1981aleph}
Richard Laver and Saharon Shelah.
\newblock The $\aleph_2$-{S}ouslin hypothesis.
\newblock {\em Transactions of the American Mathematical Society}, 264(2):411--417, 1981.

\bibitem{mekler1990universal}
Alan~H Mekler.
\newblock Universal structures in power $\aleph_1$.
\newblock {\em The Journal of symbolic logic}, 55(2):466--477, 1990.

\bibitem{scott2014logic}
Dana Scott.
\newblock Logic with denumerably long formulas and finite strings of quantifiers.
\newblock In {\em The theory of models}, pages 329--341. Elsevier, 2014.

\bibitem{shelah1990classification}
Saharon Shelah.
\newblock {\em Classification theory: and the number of non-isomorphic models}, volume~92.
\newblock Elsevier, 1990.

\bibitem{shelah1990universal}
Saharon Shelah.
\newblock Universal graphs without instances of {CH}: {Revisited}.
\newblock {\em Isr. J. Math.}, 70(1):69--81, 1990.

\bibitem{todorcevic2007lipschitz}
Stevo Todorcevic.
\newblock Lipschitz maps on trees.
\newblock {\em Journal of the Institute of Mathematics of Jussieu}, 6(3):527--556, 2007.

\bibitem{vaananen1995games}
Jouko V{\"a}{\"a}n{\"a}nen.
\newblock Games and trees in infinitary logic: A survey.
\newblock In {\em Quantifiers: Logics, Models and Computation: Volume One: Surveys}, pages 105--138. Springer, 1995.

\end{thebibliography}
\end{document}